\documentclass[a4paper, 11pt]{article}
\usepackage{pgf, tikz}
\usepackage{latexsym}
\usepackage{amsfonts, amssymb,  mathtools, enumerate, xcolor,amsthm}
\usepackage{geometry}
\usepackage{amsmath, amsfonts}
\usepackage{placeins}
\usepackage{setspace}
\usepackage{graphicx}
\usepackage{tabularx}
\numberwithin{equation}{section}
\usepackage{cite}
\usepackage{booktabs,rotating}
\usepackage{epstopdf}
\usepackage{csquotes}
\usepackage{subfig}
\usepackage{parskip}
\usepackage{manyfoot}%
\usepackage{authblk}
\usepackage{float} 
\usepackage{textcomp,booktabs}
\usepackage[colorlinks=true, linkcolor=brown, citecolor= black, filecolor=blue, urlcolor=blue]{hyperref}
\newtheorem{theorem}{Theorem}[section]
\newtheorem{lemma}{Lemma}[section]

\usepackage[toc,page,titletoc]{appendix}
\graphicspath{{figures/}}  

\DeclareCaptionFormat{cont}{#1 (cont.)#2#3\par}

\title{%
A multiscale framework integrating within-host infection kinetics with airborne transmission dynamics
}

\author{%
Andrew Omame\thanks{Department of Mathematics and Statistics, York University,
 Toronto, Ontario, Canada.} \thanks{Disease-Informed Modelling, Methods, and Systems (DIMMS) Laboratory, Department of Mathematics and Statistics, York University, Keele Campus, 4700 Keele Street, Toronto, M3J 1P3, Ontario, Canada.} \quad and \quad
Sarafa A. Iyaniwura\thanks{Vaccine and Infectious Disease Division, Fred Hutchinson Cancer Center, Seattle, WA, USA. (iyaniwura@aims.ac.za).}
}

\begin{document}
	\maketitle	

\begin{abstract}

Coupling within-host infection dynamics with population-level transmission remains a major challenge in infectious disease modeling, especially for airborne pathogens with potential to spread indoor. The frequent emergence of such diseases highlight the need for integrated frameworks that capture both individual-level infection kinetics and between-host transmission. While analytical models for each scale exist, tractable approaches that link them remain limited.
In this study, we present a novel multiscale mathematical framework that integrates within-host infection kinetics with airborne transmission dynamics. The model represents each host as a patch and couples a system of ordinary differential equations (ODEs) describing in-host infection kinetics with a diffusion-based partial differential equation (PDE) for airborne pathogen movement in enclosed spaces. These scales are linked through boundary conditions on each patch boundary, representing viral shedding and inhalation.
Using matched asymptotic analysis in the regime of intermediate diffusivity, we derived a nonlinear ODE model from the coupled ODE--PDE system that retains spatial heterogeneity through Neumann Green’s functions. We established the existence, uniqueness, and boundedness of solutions to the reduced model and analyzed within-host infection kinetics as functions of the airborne pathogen diffusion rate and host spatial configuration. In the well-mixed limit, the model recovers the classical target cell limited viral dynamics framework. 
Overall, the proposed multiscale modeling approach enables the simultaneous study of transient within-host infection dynamics and population-level disease spread, providing a tractable yet biologically grounded framework for investigating airborne disease transmission in indoor environments.
\end{abstract}
	
	
\section{Introduction}\label{sec:Intro}

Airborne infectious diseases remain a significant public health concern due to their high transmissibility and capacity to cause large-scale outbreaks. These diseases are transmitted via pathogen-laden particles that can remain suspended in the air for extended periods and be inhaled by susceptible individuals. Examples of classic airborne diseases include tuberculosis, measles, and chickenpox, while respiratory infections such as influenza, SARS-CoV-2, and adenoviruses also spread via droplets or aerosols, especially in poorly ventilated indoor environments \cite{tellier2019recognition, xie2007far, tang2006airborne}.
Indoor environments such as classrooms, hospitals, offices, and public transportation systems are particularly conducive to airborne disease transmission. In these settings, pathogen-laden aerosols can accumulate due to limited air exchange and the close proximity of individuals, leading to a higher risk of infection \cite{eames2009airborne, ai2018airborne}. The efficiency of transmission is modulated by factors such as the rate of viral shedding, aerosol properties (e.g., particle size and diffusivity), air circulation, and environmental conditions such as temperature and humidity \cite{tang2006airborne, wei2016airborne}.

Infectious disease modeling has traditionally been carried out at two distinct biological scales: within-host models that describe pathogen dynamics inside an individual, and between-host models that capture disease spread at the population level. Within-host models often employ systems of ordinary differential equations (ODEs) to characterize interactions among pathogens, host cells, and immune responses, offering insights into viral load kinetics and immune control \cite{perelson2002modelling, heitzman2022modeling, ke2021vivo, perelson2021mechanistic}. In contrast, between-host models, such as the classical susceptible–infected–recovered (SIR) framework, describe the flow of individuals through epidemiological compartments \cite{hethcote2000mathematics, kermack1927contribution}.
Increasingly, research has focused on linking these scales to build multiscale models that incorporate individual infection dynamics into population-level transmission. A growing body of work shows that within-host factors, such as viral load, immune response, and treatment, can significantly influence transmission probability, epidemic size, and pathogen evolution \cite{gilchrist2006evolution, martcheva2015coupling, mideo2008linking, alizon2009virulence}. For instance, viral load has been used as a proxy for infectiousness in both deterministic and stochastic models, and several studies have explored how host-level immunity or therapy modifies transmission parameters \cite{alizon2005emergence, wang2022multiscale}.

One of the foundational contributions in this area is by Gilchrist and Sasaki \cite{gilchrist2002modeling}, who employed an evolutionary framework in which pathogen replication within a host influences both transmission potential and host mortality, thereby linking the scales through pathogen fitness \cite{gilchrist2002modeling}. Mideo et al. \cite{mideo2008linking} further reviewed these multiscale models, highlighting how within-host dynamics drive key epidemiological traits such as virulence and transmissibility.
An alternative strategy embeds within-host ODE systems into agent-based or structured population models, treating each host as a distinct entity with internal infection dynamics. This framework enables the modeling of heterogeneity in disease progression and treatment response across individuals \cite{day2006biological, handel2010toward}.
More recent hybrid approaches incorporate dynamic transmission rates that evolve with within-host infection states. For instance, Alizon and van Baalen \cite{alizon2005emergence} proposed a model where viral load directly informs time-varying infectiousness, capturing transmission as an emergent property of within-host processes.

Despite these advances, most existing multiscale models either assume well-mixed environments or rely on simplified representations of airborne transmission. Few models explicitly account for the spatial distribution of hosts and airborne pathogens, and their impact on the local infection risk of individuals. To bridge this gap, we propose a novel mathematical framework that integrates within-host viral dynamics with a spatially resolved model of pathogen dispersion in air. In our model, hosts are represented as small patches in the environment; the spatial density of particles is modeled using a diffusion partial differential equation (PDE), while the within-host viral kinetics are described by an ODE system.
The coupling between these two scales is achieved via boundary conditions and nonlocal integrals that quantify the shedding and inhalation of viral particles. We employ matched asymptotic analysis in the limit of intermediate viral diffusivity to derive a set of nonlinear ODEs from the coupled ODE-PDE model. The resulting reduced ODE system incorporates spatial heterogeneity and host locations through a Green’s function formulation, enabling analysis of how spatial structure influences within-host infection dynamics and airborne transmission.

The rest of the article is organized as follows. In Section~\ref{sec:MathModelForm}, we present the full model formulation and its non-dimensionalization. Section~\ref{sec:AsymptoticAnalysis} describes the asymptotic reduction of the coupled ODE-PDE model to a multiscale ODE system. In Section~\ref{sec:ExistUnique}, we establish the existence, uniqueness, and boundedness of the ODE model solutions. Section~\ref{sec:ViralProperties} investigates the within-host infection kinetics using the reduced ODE system, and shows the effect of diffusion and host position on the viral kinetics. Finally, Section~\ref{sec:Discussion} discusses the significance and limitations of our approach, and outlines directions for future work.

\section{Model formulation}\label{sec:MathModelForm}

\begin{figure}[!h]
    \centering
    \includegraphics[scale=0.350]{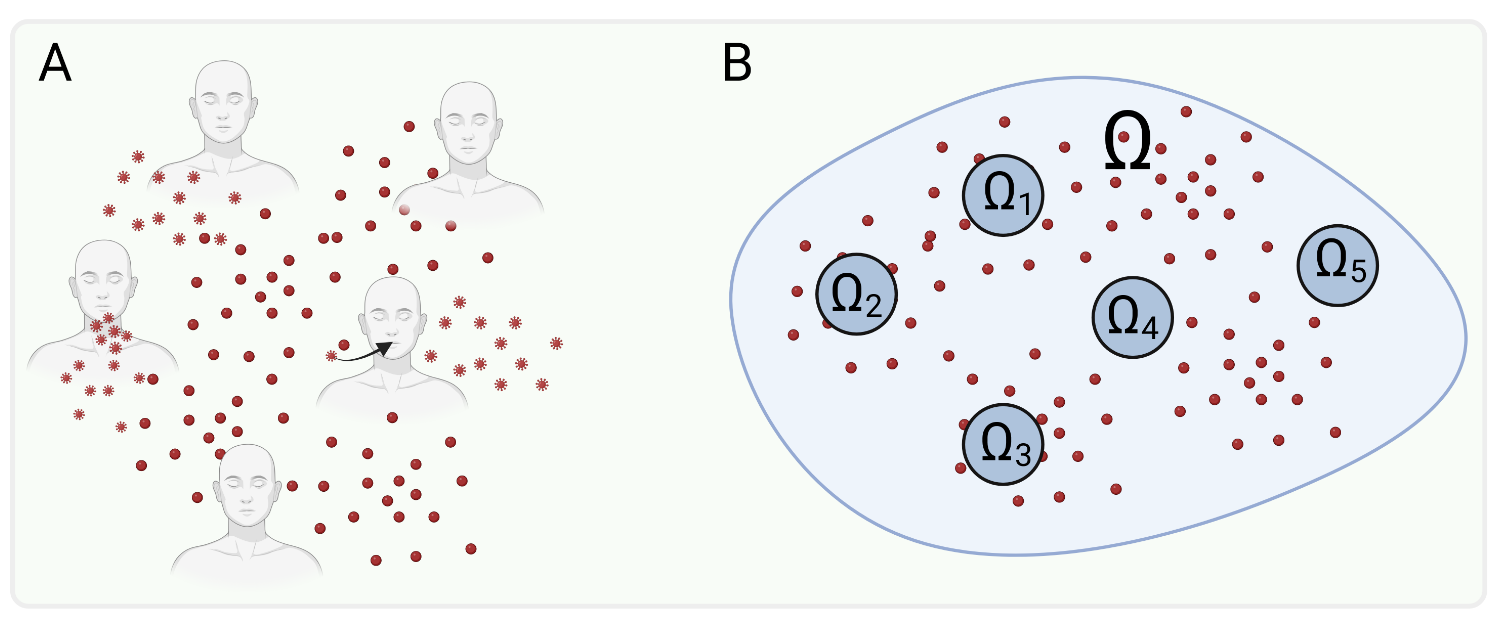}
    \caption{\textbf{Schematic illustration of our modeling framework.}  (A) Illustration of individuals in a region showing virus exhalation, diffusion, and inhalation. (B) Mathematical description  of the problem, where individuals are represented by circular patches ($\Omega_j$, for $j=1,\dots,5$), located in a bounded domain ($\Omega$). The red dots represent diffusing viral particles in both panels. Diagram is created by BioRender. }
    \label{fig:Sche_diagram_RNA}
\end{figure}

We develop a multiscale mathematical modeling framework to study the spatial spread of airborne infectious diseases among hosts in enclosed environments. The model couples the {within-host} dynamics of pathogen replication with its {between-host} transmission dynamics via airborne diffusion.
Our framework is structurally similar to reaction-diffusion models used in the study of intercellular communication \cite{gou2016asymptotic, iyaniwura2021synchronous, iyaniwura2021synchrony}, as well as to PDE-based SIR models developed for analyzing airborne disease transmission across populations \cite{david2020novel, david2022effect}. We assume that hosts are localized and spatially distributed in a confined domain (e.g., hospital wards, prison cells, office spaces, or classrooms), where infected individuals release airborne pathogens, such as viruses, into the surrounding environment.
The spatial dynamics of the airborne pathogen are modeled by a linear diffusion PDE, while the within-host infection kinetics are captured using a system of ODEs (see Figure~\ref{fig:Sche_diagram_RNA}). These two scales are coupled through a nonlocal term that integrates the pathogen concentration in the vicinity of each host.
Although we focus on viral infections for the sake of illustration, the modeling framework is general and applicable to a wide range of airborne pathogens. Throughout the article, we use the terms \textit{host} and \textit{individual} interchangeably.

Let $\Omega \subset \mathbb{R}^2$ be a bounded domain representing an enclosed region in which $m$ individuals are located, and let each individual be represented by $\Omega_j \subset \Omega$, for $j = 1, \dots, m$, as illustrated in Figure~\ref{fig:Sche_diagram_RNA}B. The surrounding space, $\Omega \setminus \bigcup_{j=1}^{m} \Omega_j$, represents the region of air where viral particles secreted by infected individuals can diffuse. 
We denote the density of virus particles in the air at position $\pmb{X} \in \Omega \setminus \bigcup_{j=1}^{m} \Omega_j$ and time $\tau$ by $r(\pmb{X}, \tau)$. This density evolves over time due to diffusion and is influenced by viral shedding from infected hosts (see Figure~\ref{fig:Sche_diagram_RNA}A). The dynamics of $r(\pmb{X}, \tau)$ is governed by the following PDE:
\begin{subequations}\label{eqn:Dim_BulkModel}
\begin{equation}\label{eqn:Dim_Bulk}
    \begin{split}
 \partial_{\tau} r &= D_r \,\Delta_X r - k_r\, r, \quad \tau > 0 , \quad X \in \Omega \setminus \cup_{j=1}^{m} \Omega_j \,; \\[2ex]
 \partial_{n_X} r = 0, & \quad \pmb{X} \in  \partial \Omega,  \qquad  D_r\, \partial_{n_X} r = - \gamma_j \,u_j, \quad \pmb{X} \in \partial \Omega_j, \qquad  j=1,\dots,m.
    \end{split}
\end{equation}
where $D_r > 0$ and $k_r > 0$ denote the dimensional diffusion and degradation rates of viral particles in the air, respectively. The variable $u_j \equiv u_j(\tau)$ represents the viral concentration within the $j^{\text{th}}$ host at time $\tau$, and $\gamma_j$ is the rate at which virus particles are exhaled into the air by the $j^{\text{th}}$ host.
The operator $\partial_{n_{\pmb{X}}}$ denotes the outward normal derivative on the boundary $\partial \Omega_j$ of the $j^{\text{th}}$ host, pointing into the surrounding environment. The Robin boundary condition imposed in Eq.~\eqref{eqn:Dim_Bulk} indicates that the rate at which virus particles are expelled into the air is proportional to the within-host viral concentration. Consequently, if the $j^{\text{th}}$ host is uninfected (i.e., $u_j = 0$), then no virus is exhaled into the air. In contrast, an infected individual continuously releases viral particles at a rate proportional to their viral load, modulated by $\gamma_j$, until the end of their infection.
Once expelled, these viral particles diffuse throughout the domain $\Omega \setminus \bigcup_{j=1}^{m} \Omega_j$ and degrade over time, as described by the PDE in Equation~\eqref{eqn:Dim_Bulk} and illustrated in Figure~\ref{fig:Sche_diagram_RNA}A.

To capture the infection kinetics within an infected host, we adopt a compartmental framework inspired by previous viral dynamics models \cite{ke2021vivo, perelson2021mechanistic, iyaniwura2024kinetics}. The model divides host cells into three key compartments: target cells ($\mathcal{T}$), representing uninfected and susceptible cells; eclipse phase cells ($\mathcal{E}$), which are infected but not yet producing virus; and productively infected cells ($\mathcal{I}$), which actively produce and release viral particles. The eclipse phase accounts for the intracellular delay between initial infection and the onset of viral production \cite{kakizoe2015method, vaidya2021modeling}, analogous to the exposed compartment in SEIR models commonly used in population-level epidemiology.
We also track the concentration of free virus within the host, denoted by $u_j$. The within-host viral load in each individual ($u_j$) is coupled to the airborne virus dynamics through boundary conditions (Eq.~\eqref{eqn:Dim_Bulk}), thereby linking individual infection kinetics to population-level transmission.
This framework assumes a fixed pool of target cells at the onset of infection, with no replenishment during the infection course, consistent with previous studies \cite{baccam2006kinetics, best2018mathematical, de1998target}. The complete set of within-host model equations is given by:
\begin{equation}\label{eqn:basic_viral_eclipse_dyn}
\begin{split}
        \frac{\text{d}\mathcal{T}_j}{\text{d}\tau} &= -b_{1 j} \mathcal{T}_j \int_{\partial \Omega_j} \Big{(} r/r_c \Big{)} \, \text{d}S_{X} - b_{2 j} \mathcal{T}_j \, u_j, \\[1ex]
\frac{\text{d}\mathcal{E}_{j}}{\text{d}\tau} & = b_{1j} \mathcal{T}_j \int_{\partial \Omega_j} \Big{(} r/r_c \Big{)} \, \text{d}S_{X} + b_{2j} \mathcal{T}_j \, u_j - \alpha_{j} \mathcal{E}_{j},\\[1ex]
\frac{\text{d}\mathcal{I}_{j}}{\text{d}\tau} & = \alpha_{j} \, \mathcal{E}_{j} - d_j \, \mathcal{I}_{j},   \\[1ex]
\frac{\text{d} u_{j}}{\text{d}\tau} & =  \rho_j \, \mathcal{I}_{j} - \phi_j \, u_j, \\
\end{split}    
\end{equation}
\end{subequations}
where \( b_{1j} \) and \( b_{2j} \) are the dimensional infection rates of the diffusing and in vivo viruses, respectively, \( d_j \) is the dimensional death rate of infected cells, \( \alpha_j \) is the dimensional rate at which infected cells transition from the eclipse phase to the virus-producing phase, \( \rho_j \) is the dimensional virus production rate by infected cells, \( \phi_j \) is the dimensional in vivo virus clearance rate, and \( r_c \) is a typical value for the density of viruses in the air.
Recall that \( r \) is the density of diffusing viral particles, whose dynamics are governed by the PDE in \eqref{eqn:Dim_Bulk}. The integrals in \eqref{eqn:basic_viral_eclipse_dyn} account for the density of viral particles around the \( j^{\text{th}} \) host, which determines the initiation of infection in a susceptible individual.
In other words, our model assumes that a susceptible individual becomes infected based on the density of viral particles around them. The susceptible cells in each host are infected either by the diffusing viruses inhaled by the host or by the viruses produced within the host. These two infection routes are represented by the two terms in the equation for \( \mathcal{T}_j \) in \eqref{eqn:basic_viral_eclipse_dyn}.

Next, we non-dimensionalize the multiscale model in \eqref{eqn:Dim_BulkModel}. First, we define the dimensions of the variables in the model as follows:
\begin{equation}
    \begin{split}
[r] &=  \frac{\text{copies/mL}}{\text{length}^2}, \qquad [D_r] = \frac{\text{length}^2}{\text{time}}, \qquad [r_c] = \text{copies/mL}, \qquad [\tau] = \text{time}, \quad [X] = \text{length}, \\[1ex]
[b_{1j}] &=\frac{\text{length}}{\text{time}}, \qquad [b_{2j}] = \frac{1}{\text{time} \times \text{copies/mL}}, \qquad 
[\mathcal{T}_j] = [\mathcal{E}_{j}] = [\mathcal{I}_{j}] = \text{cells}, \qquad [u_j] = \text{copies/mL}, \\
[\alpha_j] &= [d_j] = [k_r] = [\phi_j] =  \frac{1}{\text{time}}, \quad [\rho_j] = \frac{\text{copies/mL}}{\text{time} \times \text{cells}}, \quad
[\gamma_{j}] =  \frac{1}{\text{length} \times \text{time} }, 
    \end{split}
\end{equation}
for $ j=1,\dots,m$, where $[y]$ represents the dimension of $y$.

Throughout this article, we represent the region of interest (our domain) with a disk of radius $L$ and the hosts as smaller disks of common radius $R$. We assume that the hosts are relatively small, compared to the size of the domain and introduce a scaling parameter $\varepsilon = R/L \ll 1$. We define the dimensionless variables as 
\begin{equation}\label{eqn:dimless_var}
\begin{split}
   \mathcal{V} &= \frac{L^2}{r_c} r, \qquad T_j = \frac{\mathcal{T}_j}{N_j},  \qquad E_j = \frac{\mathcal{E}_j}{N_j}, \qquad  I_{j} = \frac{\mathcal{I}_{j}}{N_j}, \qquad v_j = \frac{u_j}{r_c}, \qquad \pmb{x} = \frac{\pmb{X}}{L}, \qquad t = k_r \tau,
\end{split}    
\end{equation}
where $N_j = \mathcal{T}_j +  \mathcal{E}_j + \mathcal{I}_j$ is the total number of cells (relating to this disease) in the $j^{\text{th}}$ host, and $r_c$ is the typical values for $r$. Based on the definitions in \eqref{eqn:dimless_var}, $T_j$, $E_{ j}$, and $I_j$, are the proportion of target cells, cells in the eclipse phase, and productively infected cells, respectively, in the $j^{\text{th}}$ host.

We use the dimensionless variables in \eqref{eqn:dimless_var} to non-dimensionalize the multiscale model \eqref{eqn:Dim_BulkModel}. In dimensionless form, the dynamics of the diffusing viral particles is given by
\begin{subequations}\label{eqn:Dimless_Coupled_EclipDyn}
\begin{equation}
\begin{split}
    \partial_t \mathcal{V} &= D  \, \Delta \mathcal{V} - \mathcal{V}, \quad t > 0 , \quad \pmb{x} \in \Omega \setminus \cup_{j=1}^{m} \Omega_{\varepsilon j}; \label{eqn:Dimless_Coupled_EclipDynA}\\[1ex]
\partial_{n} \mathcal{V} = 0, \quad \pmb{x} \in & \partial \Omega\, ;  \qquad D\, \partial_{n} \mathcal{V} = - \xi_j \,v_j/|\partial \Omega_{\varepsilon j}|, \quad \pmb{x} \in \partial \Omega_{\varepsilon j}, \qquad  j=1,\dots,m. 
\end{split}
\end{equation}
This dynamics is coupled to the dimensionless dynamics of the $j^{\text{th}}$ host given by 
\begin{equation}\label{eqn:Dimless_Coupled_EclipDynB}
\begin{split}
\frac{\text{d} T_j}{\text{d}t} &= - \frac{\beta_{1j}\, T_j}{|\partial\Omega_{\varepsilon j}|}  \int_{\partial \Omega_{\varepsilon j}} \mathcal{V} \; \text{d}S_{\pmb{x}} - \beta_{2j} T_j \,v_j,\\[1ex] 
\frac{\text{d} E_{j}}{\text{d}t} & = \frac{\beta_{1j} \, T_j}{|\partial\Omega_{\varepsilon j}|} \int_{\partial \Omega_{\varepsilon j}} \mathcal{V} \;  \text{d}S_{\pmb{x}} + \beta_{2j}  T_j \,v_j - k_j E_j,\\[1ex] 
\frac{\text{d} I_{ j}}{\text{d}t} & = k_j E_j - \delta_j I_j,\\[1ex] 
\frac{\text{d} {v}_j}{\text{d} t} & = p_j {I}_j - c_j {v}_j, 
\end{split}
\end{equation}
\end{subequations}
where $\Omega_{\varepsilon j}$ represents the $j^{\text{th}}$ host located at $\pmb{x}_j$ with boundary $\partial \Omega_{\varepsilon j}$ and $|\partial \Omega_{\varepsilon j}| = 2\pi\varepsilon$ is its perimeter. Here, $D = D_r/(k_r L^2)$ is the dimensionless diffusion rate of the viruses in the air, $k_j = \alpha_j/k_r$ is the dimensionless rate at which cells in eclipse phase become virus producing cells, $\delta_j = d_j/k_r$ is the dimensionless death rate of infected cells, $p_j = \rho_j N_j/ (k_r\,r_c)$ is the dimensionless virus production rate by infected cells,  and $c_j= \phi_j/k_r$ is the dimensionless clearance rate of viral particles in the body of the $j^{\text{th}}$ individual.

The dimensionless infectious rate of the diffusing viral particles ($\beta_{1j}$) and within-host viruses ($\beta_{2j}$) are defined as $\beta_{1j}/|\partial \Omega_{\varepsilon j}| = b_{1 j}/(k_r\, L)$ and $\beta_{2j} = b_{2 j} r_c/k_r$.
In addition, we define the dimensionless rate of expelling viruses by the  $j^{\text{th}}$ infected host as $\xi_j/|\partial \Omega_{\varepsilon j}| = (\gamma_j \,L)/k_r$. We have scaled the dimensionless infectious rate of the diffusing viral particles ($\beta_{1j}$) and the virus secretion rate ($\xi_j$) by the perimeter of the host ($|\partial \Omega_{\varepsilon j}|$) to 
account for the surfaces through which  diffusing viral particles can enter the  $j^{\text{th}}$ host, and where viruses are secreted by the  host.
 In this non-dimensionalization, we assume that $b_{1 j}/(k_r\, L)$ and $\gamma_j \,(L\, r_c)/(k_r r_c)$ are both $\mathcal{O}(1/\varepsilon)$, which ensures that $\beta_{1j}$ and $\xi_j$ are $\mathcal{O}(1)$. This scaling is necessary to properly account for disease transmission and virus secretion due to the size of the disks used to  represent hosts, relative to the size of the domain.

\section{Asymptotic analysis}\label{sec:AsymptoticAnalysis}

We use singular perturbation theory to analyze the dimensionless multiscale model in Eq.~\eqref{eqn:Dimless_Coupled_EclipDyn} in the limit of intermediate diffusing viral particles in the air, i.e., $D = \mathcal{O}(\mu^{-1})$, where $\mu = -1/\log(\varepsilon)$ with $\varepsilon \ll 1$. From the multiscale PDE model in Eq.~\eqref{eqn:Dimless_Coupled_EclipDyn}, we aim to derive a nonlinear system of ODEs that incorporates the diffusion of viral particles and an $\mathcal{O}(\mu)$ corrections term that incorporates the effect of host location into the ODE system. In this approach, we re-scale the diffusion rate of viral particles as follows
\begin{equation}\label{eqn:rescale_D}
    D = \frac{D_0}{\mu}, \quad \text{where} \quad D_0 = \mathcal{O}(1) \quad \text{and} \quad  \mu = -1/\log(\varepsilon) \quad \text{with} \quad \varepsilon \ll 1.
\end{equation}
Assuming that the $m$ hosts, $\Omega_{\varepsilon 1}, \dots, \Omega_{\varepsilon m}$, under consideration are located at the  $\pmb{x}_1, \dots, \pmb{x}_m$, respectively.
We define $V$ as the average density of viral particles in the air, $\Omega \setminus \Omega_h$, where $\Omega_h = \cup_{j=1}^{m} \Omega_{\varepsilon j}$ (bulk region), which is given by
\begin{equation}\label{eqn:AveBulk_Pathogen}
 V = \frac{1}{|\Omega \setminus \Omega_h|} \int_{\Omega \setminus \Omega_h} \mathcal{V }\,\, \text{d} \pmb{x}.
\end{equation}
Upon multiplying the PDE in \eqref{eqn:Dimless_Coupled_EclipDynA} by $1/|\Omega \setminus \Omega_h|$ and integrating over $\Omega \setminus \Omega_h$, we have
\begin{equation}\label{eqn:AveBulk_Pathogen_2}
    \partial_t V + V = \frac{1}{|\Omega \setminus \Omega_h|} \int_{\Omega \setminus \Omega_h} D \Delta V \,\, \text{d} \pmb{x}.
\end{equation}
Using the divergence theorem on the integral in the R.H.S of \eqref{eqn:AveBulk_Pathogen_2} and the boundary conditions in \eqref{eqn:Dimless_Coupled_EclipDynA}, \eqref{eqn:AveBulk_Pathogen_2} becomes
\begin{equation}\label{eqn:AveBulk_Pathogen_3}
    \partial_t V + V = \frac{1}{|\Omega \setminus \Omega_h|} \sum_{j=1}^m \xi_j\, v_j.
\end{equation}
We have derived an equation for the average number of viral particles in the bulk region. Since $|\Omega_h|$ is the union of all hosts and $|\Omega_{\varepsilon j}| = \pi \varepsilon^2$ is the surface area of the $j^{\text{th}}$ host, then the combined surface area  of all hosts is given by $|\Omega_h| = m \pi \varepsilon^2$. Based on this derivation, we deduce that $|\Omega \setminus \Omega_h| = |\Omega| - \mathcal{O}(\varepsilon^2)$, which implies that $|\Omega \setminus \Omega_h| \to |\Omega|$ as $\varepsilon \to 0$ in \eqref{eqn:AveBulk_Pathogen_3}. We shall use this derivation subsequently in our analysis.

Next, substituting \eqref{eqn:rescale_D} into \eqref{eqn:Dimless_Coupled_EclipDynA}, we obtain
\begin{equation}
\begin{split}\label{eqn:Bulk_D01}
        \partial_t \mathcal{V} &= \frac{D_0}{\mu}  \, \Delta \mathcal{V} - \mathcal{V}, \quad t > 0 , \quad \pmb{x} \in \Omega \setminus \cup_{j=1}^{m} \Omega_{\varepsilon j};  \\
\partial_{n} \mathcal{V} = 0, \quad \pmb{x} \in & \partial \Omega\, ;  \qquad  \frac{D_0}{\mu}\, \partial_{n} \mathcal{V} = - \frac{\xi_j v_j}{2\pi \varepsilon }, \quad \pmb{x} \in \partial \Omega_{\varepsilon j}, \qquad  j=1,\dots,m .
\end{split}
 \end{equation}
 Observed that we have used $|\partial \Omega_{\varepsilon j}| = 2\pi \varepsilon$,  the perimeter of each host in this derivation.
Since each infected host expels viral particles, we expect the density of these particles around the hosts to be relative higher than in the remaining part of the domain. As a result, we consider an inner region near each host and introduce the inner variables $\pmb{y} = \varepsilon^{-1}|\pmb{x} - \pmb{x}_j|$ and $V(\pmb{x}) = Q(\varepsilon \pmb{y} + \pmb{x}_j)$, with $|\pmb{y}| = \sigma$, for $j=1,\dots,m$. Upon writing the PDE in \eqref{eqn:Bulk_D01} in terms of the inner variables, close to the $j^{\text{th}}$ host,  as $\varepsilon \to 0$, we have
\begin{equation}
\begin{split}\label{eqn:Bulk_D0}
     \Delta_{\sigma} \, Q_j &= 0, \quad t > 0 , \quad \sigma >1;  \\
 2\pi  \frac{D_0}{\mu}\, \partial_{\sigma} Q_j &= - \xi_j v_j, \quad \text{on} \quad \sigma=1, \qquad  j=1,\dots,m,
\end{split}
 \end{equation}
where $\Delta_{\sigma}$ is the radially symmetric Laplace operator with respect to $\sigma$ in 2-D. We consider the following asymptotic expansion for $Q_j(\sigma,t)$ in the inner region near the $j^{\text{th}}$ host
\begin{equation}\label{eqn:Asyp_Expand}
Q_j(\sigma,t) = Q_{0 j}(\sigma,t) + \frac{\mu}{D_0} Q_{1 j}(\sigma,t) + \dots
\end{equation}
We substitute \eqref{eqn:Asyp_Expand} into \eqref{eqn:Bulk_D0} and collect terms in powers of $\mu$. At leading-order, we have
\begin{equation}\label{eqn:Lead_Order}
\begin{split}
 \Delta_{\sigma} \, Q_{0 j} &= 0, \quad t > 0 , \quad \sigma >1;  \qquad \partial_{\sigma} Q_{0 j} = 0, \quad \text{on} \quad \sigma=1, \qquad  j=1,\dots,m,
\end{split}
\end{equation}
We observe that any constant or function of time satisfies \eqref{eqn:Lead_Order}. Therefore, we set $Q_{0j} \equiv Q_{0j}(t)$. At the next-order, we have
\begin{equation}\label{eqn:Second_Order}
\begin{split}
 \Delta_{\sigma} \, Q_{1j} &= 0, \quad t > 0 , \quad \sigma >1;  \\[1ex]
 2\pi  \, \partial_{\sigma} Q_{1j} &= - \xi_j v_j, \quad \text{on} \quad \sigma=1, \qquad  j=1,\dots,m,
\end{split}
\end{equation}
which has the solution
\begin{equation}\label{eqn:Second_Order_Sol}
\begin{split}
  Q_{1j} &= \frac{-\xi_j v_j}{2\pi}\log(\sigma) + \Upsilon_j, \qquad  j=1,\dots,m,
\end{split}
\end{equation}
where $\Upsilon_j$ is  a constant to be determined. In terms of the outer variables, we have
\begin{equation}\label{eqn:Second_Order_Sol}
\begin{split}
  Q_{1j} &= \frac{-\xi_j v_j}{2\pi} \left( \ln |\pmb{x} - \pmb{x}_j| + \frac{1}{\mu} \right) + \Upsilon_j, \qquad  j=1,\dots,m,
\end{split}
\end{equation}
Upon substituting $Q_{0 j}$ and $Q_{1j}$ into the asymptotic expansion \eqref{eqn:Asyp_Expand}, we have 
\begin{equation}\label{eqn:Asyp_Expand_2}
Q_j = \left( Q_{0 j} - \frac{\xi_j v_j}{2\pi D_0} \right) + \frac{\mu}{D_0} \left( -\frac{\xi_j v_j}{2\pi} \ln |\pmb{x} - \pmb{x}_j| + \Upsilon_j  \right) + \dots
\end{equation}

Now, we consider the outer problem at an $\mathcal{O}(1)$ distance from each host
\begin{equation}
\begin{split}\label{eqn:Outer}
        \partial_t \mathcal{V} &= \frac{D_0}{\mu}  \, \Delta \mathcal{V} - \mathcal{V}, \quad t > 0 , \quad \pmb{x} \in \Omega \setminus \{\pmb{x}_1,\dots,\pmb{x}_m  \};  \\
&\partial_{n} \mathcal{V} = 0, \quad \pmb{x} \in  \partial \Omega,  \qquad 
\end{split}
 \end{equation}
 where the positions of the hosts are represented by $\pmb{x}_1,\dots,\pmb{x}_m$. In the outer region, we expand the outer solution $\mathcal{V}$ as
 \begin{equation}\label{eqn:Outer_expand}
     \mathcal{V} = \mathcal{V}_0 + \frac{\mu}{D_0} \mathcal{V}_1 + \dots
 \end{equation}
 Upon substituting the asymptotic expansion \eqref{eqn:Outer_expand} into \eqref{eqn:Outer} and collecting terms in powers of $\mu$, at leading-order, we have
 \begin{equation}
\begin{split}\label{eqn:Lead_Outer}
         \Delta \mathcal{V}_0 = 0, \quad t > 0 , \quad \pmb{x} \in \Omega \setminus \{\pmb{x}_1,\dots,\pmb{x}_m  \};  \qquad
\partial_{n} \mathcal{V}_0 = 0, \quad \pmb{x} \in  \partial \Omega. 
\end{split}
 \end{equation}
 Similar to the leading-order inner problem \eqref{eqn:Lead_Order}, any constant or function of time satisfies \eqref{eqn:Lead_Outer}. Therefore, we set $\mathcal{V}_0 \equiv \mathcal{V}_0(t)$. At $\mathcal{O}(\mu)$, we have
 \begin{equation}
\begin{split}\label{eqn:Outer_Order_Nu}
    \Delta \mathcal{V}_1 =  \mathcal{V}_0 + \mathcal{V}_{0 t}, \quad t > 0 , \quad \pmb{x} \in \Omega \setminus \{\pmb{x}_1,\dots,\pmb{x}_m  \};  \qquad
\partial_{n} \mathcal{V}_1 = 0, \quad \pmb{x} \in  \partial \Omega.
\end{split}
 \end{equation}
Matching the inner expansion \eqref{eqn:Asyp_Expand_2} to the outer expansion \eqref{eqn:Outer_expand}, we obtained the matching conditions
 \begin{subequations}
 \begin{align}
     \mathcal{V}_0 & \sim Q_{0 j} - \frac{\xi_j \, v_j}{2 \pi D_0}, \label{eqn:Singular_1}\\[1ex]
     \mathcal{V}_1 & \sim  - \frac{\xi_j \, v_j}{2 \pi } \log |\pmb{x} - \pmb{x}_j| + \Upsilon_j, \quad \text{as} \quad \pmb{x} \to \pmb{x}_j \quad \text{for} \quad j=1,\dots,m. \label{eqn:Singular_2}
 \end{align}
 \end{subequations}
Considering the $\mathcal{O}(\mu)$ outer problem \eqref{eqn:Outer_Order_Nu} together with the singularity behavior of $\mathcal{V}_1$ in \eqref{eqn:Singular_2}, we construct the complete outer problem for $\mathcal{V}_1$ given by
 \begin{equation}\label{eqn:OuterComplete}
\begin{split}
    \Delta \mathcal{V}_1 &=  \mathcal{V}_0 + \mathcal{V}_{0 t}, \quad t > 0 , \quad \pmb{x} \in \Omega; \qquad
\partial_{n} \mathcal{V}_1 = 0, \quad \pmb{x} \in  \partial \Omega; \\ 
 & \mathcal{V}_1  \sim  - \frac{\xi_j \, v_j}{2 \pi} \log |\pmb{x} - \pmb{x}_j| + \Upsilon_j, \quad \text{as} \quad \pmb{x} \to \pmb{x}_j.
\end{split}
 \end{equation}
 Upon writing the singularity behavior of $\mathcal{V}_1$ in \eqref{eqn:OuterComplete} as a delta function in the PDE for $\mathcal{V}_1$, we have
\begin{equation}\label{eqn:OuterComplete_2}
\begin{split}
    \Delta \mathcal{V}_1 &=  \mathcal{V}_0 + \mathcal{V}_{0 t} - \sum_{i=1}^{m} \xi_i \,v_i \, \delta( \pmb{x} - \pmb{x}_i), \quad t > 0 , \quad \pmb{x} \in \Omega; \qquad
\partial_{n} \mathcal{V}_1 = 0, \quad \pmb{x} \in  \partial \Omega.
\end{split}
 \end{equation} 
The next step in this analysis is to derive a solvability condition for the complete outer problem \eqref{eqn:OuterComplete_2}. We integrate the PDE in \eqref{eqn:OuterComplete_2} over the domain $\Omega$ to obtain
 \begin{equation}\label{eqn:OuterComplete_Integral}
\begin{split}
 \int_{\Omega}   \Delta \mathcal{V}_1 \,\, \text{d}\pmb{x} &=  (\mathcal{V}_0 + \mathcal{V}_{0 t})|\Omega| - \sum_{i=1}^{m} \xi_i \,v_i .
\end{split}
 \end{equation}
 Using the divergence theorem together with the boundary condition that  $\partial_n \mathcal{V}_1 = 0$ on $\partial \Omega$, we obtain the singularity behaviour for \eqref{eqn:OuterComplete_2}, which is given by
 \begin{equation}\label{eqn:Bulk_ODE_1}
\begin{split}
    \mathcal{V}_{0 t} &= - \mathcal{V}_0 + \frac{1}{|\Omega|} \sum_{i=1}^{m} \xi_i v_i.
\end{split}
 \end{equation}
 Observe that \eqref{eqn:Bulk_ODE_1} is the same as the ODE obtained for the average density of viral particles in \eqref{eqn:AveBulk_Pathogen_3}. Imposing the zero average constraint on $\mathcal{V}_1$, i.e., $\int_{\Omega} \mathcal{V}_1 \; \text{d}\pmb{x} = 0 $ ensures that $\mathcal{V}_0$ is the spatial average of $\mathcal{V}$ up to the $\mathcal{O}(\mu)$ term.
 The summation in the ODE \eqref{eqn:Bulk_ODE_1} shows the averaging of the contribution of virus from each host into the air.
 
 Next, we write the solution of the $\mathcal{O}(1)$ problem \eqref{eqn:OuterComplete_2}, which satisfies $\int_{\Omega} \mathcal{V}_1 \; \text{d}\pmb{x} = 0 $, as
 \begin{equation}\label{eqn:Outer_Sol}
     \mathcal{V}_1  = \sum_{i=1}^{m} \xi_i \, v_i \, G(\pmb{x}, \pmb{x}_i).
 \end{equation}
 In \eqref{eqn:Outer_Sol},  $G(\pmb{x}, \pmb{x}_j)$ is the Neumann Green's function satisfying
 \begin{equation}\label{eqn:Greens_func}
\begin{split}
    \Delta G &= \frac{1}{|\Omega|} - \delta(\pmb{x} - \pmb{x}_j), \quad \pmb{x} \in \Omega; \qquad \partial G = 0, \quad \pmb{x} \in \partial \Omega;\\
     G(\pmb{x}, \pmb{x}_j) & \sim -\frac{1}{2 \pi} \ln |\pmb{x} - \pmb{x}_j| + R_j, \quad \text{as} \quad \pmb{x} \to \pmb{x}_j\, \quad \text{and} \quad \int_{\Omega} G \,\, \text{d}\pmb{x} = 0,
\end{split}     
 \end{equation}
 where $R_j \equiv R(\pmb{x}_j)$ is the regular part of $G(\pmb{x}, \pmb{x}_j)$ as $\pmb{x} \to \pmb{x}_j$. 
 The solution of this PDE in a disk region is well established \cite{kolokolnikov2005optimizing, david2020novel, iyaniwura2021synchrony} with the regular part given by
\begin{equation}\label{Eq:GreensFunction_Reg}
    \begin{split}
G(\pmb{x}; \pmb{x}_j) &= -\frac{1}{2\pi} \log |{\pmb{x}} - \pmb{x}_j| 
                       - \frac{1}{4\pi} \log\left(|{\pmb{x}}|^2 |\pmb{x}_j|^2 + 1 - 2 {\pmb{x}} \cdot \pmb{x}_j \right) 
                       + \frac{(|{\pmb{x}}|^2 + |\pmb{x}_j|^2)}{4\pi} - \frac{3}{8\pi},
    \end{split}
\end{equation}
and the singular part given by
\begin{equation}\label{Eq:GreensFunction_Sing}
    \begin{split}
R(\pmb{x}_j) = -\frac{1}{2\pi} \log(1 - |\pmb{x}_j|^2) + \frac{|\pmb{x}_j|^2}{2\pi} - \frac{3}{8\pi}.    \end{split}
\end{equation}
Expanding the outer solution \eqref{eqn:Outer_Sol} as $\pmb{x} \to \pmb{x}_j$ and using the singularity behaviour of the Neumann Green's function given in \eqref{eqn:Greens_func}, we obtain
 \begin{equation}\label{eqn:P1_singular}
     \mathcal{V}_1 \sim -\frac{\xi_j v_j}{2 \pi} \ln |\pmb{x} - \pmb{x}_j| + \xi_j v_j R_j +
     \sum_{i\neq j}^{m} \xi_i \, v_i \, G(\pmb{x}_j, \pmb{x}_i) \quad \text{as} \quad \pmb{x} \to \pmb{x}_j\,.
 \end{equation}
Upon matching the singularity behaviour of $\mathcal{V}_1$ given in \eqref{eqn:OuterComplete} together with \eqref{eqn:P1_singular}, we deduced that 
\begin{equation}\label{eqn:k}
   \Upsilon_j \sim \xi_j v_j R_j + 
     \sum_{i\neq j}^{m} \xi_i \, v_i \, G(\pmb{x}_j, \pmb{x}_i) \quad \text{as} \quad \pmb{x} \to \pmb{x}_j\,.
\end{equation}
 We substitute $\Upsilon_j$ into the inner expansion \eqref{eqn:Asyp_Expand_2} to obtain a two term asymptotic expansion for the inner solution near the $j^{\text{th}}$ host,  which is given by
\begin{equation}\label{eqn:Inner_Sol}
Q_j = \left( Q_{0 j} - \frac{\xi_j v_j}{2\pi D_0} \right) + \frac{\mu}{D_0} \left( -\frac{\xi_j v_j}{2\pi} \ln |\pmb{x} - \pmb{x}_j| + \xi_i v_i R_j + 
     \sum_{i\neq j}^{m} \xi_i \, v_i \, G(\pmb{x}_j, \pmb{x}_i)  \right) + \dots, \qquad j=1,\dots,m.
\end{equation}
In matrix form, we have
 \begin{equation}\label{eqn:Inner_Sol_matrix}
Q_j =  Q_{0 j}  + \frac{\mu}{D_0} \left( -\frac{\xi_j v_j}{2\pi} \ln (\sigma) + (\mathcal{G} \Psi)_j  \right) + \dots, \qquad j=1,\dots,m.
\end{equation}
Note that we have used $\sigma = |\pmb{x} - \pmb{x}_j|/\varepsilon$ and $\mu = -1/\log(\varepsilon)$ in deriving  \eqref{eqn:Inner_Sol_matrix}. Here, $(\mathcal{G} \Psi)_j$ is the $j^{\text{th}}$ entry of the matrix-vector product $\mathcal{G} \Psi$, with $\Psi = (\xi_1 v_1, \dots, \xi_m v_m)$ and $\mathcal{G}$ is the Neumann Green's matrix define by
 \begin{equation}\label{eqn:G_matrix}
     \mathcal{G} =  
     \begin{pmatrix}
      R_1 & G_{1,2} & G_{1,3} & \dots & \dots & G_{1,m} \\
      G_{2,1} & R_2 & G_{2,3} & \dots & \dots & G_{2,m}  \\
      \vdots & \ddots & \ddots & \ddots & \dots & \vdots \\
      G_{m-2,1} &  \dots & G_{m-2,m-3} & R_{m-2} & G_{m-2,m-1} & G_{m-1,m} \\
      G_{m-1,1} &  \dots & \dots & G_{m-1,m-2} & R_{m-1} & G_{m-1,m} \\
      G_{m,1} & G_{m,2} & \dots & \dots &  G_{m,m-1} & R_{m}
     \end{pmatrix},
 \end{equation}
where $G(\pmb{x}_j;\pmb{x}_i)$ is the Neumann Green's function satisfying \eqref{eqn:Greens_func} and $R_j$ is its regular part at $\pmb{x}_j$. 

Since $\sigma=1$ on the boundary of the disk representing the  $j^{\text{th}}$ host and 
 $ Q_{0 j} = \mathcal{V}_0 + \frac{\xi_j \, v_j}{2 \pi D_0}$ from \eqref{eqn:Singular_1},  \eqref{eqn:Inner_Sol_matrix} reduces to
\begin{equation}\label{eqn:Inner_Sol_matrix_2}
Q_j = \left( \mathcal{V}_0 + \frac{\xi_j v_j}{2\pi D_0} \right) + \frac{\mu}{D_0}  \Big{(} \mathcal{G} \Psi \Big{)}_j  + \dots, \qquad j=1,\dots,m.
\end{equation}
The inner solution $Q_j$ in \eqref{eqn:Inner_Sol_matrix_2} gives an expression for the density of viral particles near the $j^{\text{th}}$ host. This implies that $\mathcal{V} \sim Q_j$ as $\pmb{x} \to \pmb{x}_j$. 

Now, integrating the asymptotic expansion for $\mathcal{V }(t)$ in \eqref{eqn:Outer_expand} over the bulk region $|\Omega \setminus \Omega_h|$ as $\varepsilon \to 0$ and using $\int_{\Omega} \mathcal{V}_1 \; \text{d}\pmb{x} = 0$, we obtain $V(t) = \mathcal{V}_0(t)$, 
which shows that $\mathcal{V}_0(t)$ is the spatial average of the virus in the bulk region far away from the cells.
Based on this, from \eqref{eqn:AveBulk_Pathogen_3}, in the limit $\varepsilon \to 0$, we have the following ODE for the average virus in the air
\begin{subequations}\label{eqn:ReducedODE_EclipDyn_IVP}
\begin{equation}\label{eqn:ReducedODE_EclipDyn_Virus}
    \frac{\text{d}V(t)}{\text{d}t} = \frac{1}{|\Omega|} \sum_{i=1}^{m} \xi_i\, v_i(t) - V(t), 
\end{equation}
which is equivalent to \eqref{eqn:Bulk_ODE_1}. Since $\mathcal{V} \sim Q_j$ as $\pmb{x} \to \pmb{x}_j$, substituting $Q_j$ as given in \eqref{eqn:Inner_Sol_matrix_2} into the ODE system for the infection kinetics of the $j^{\text{th}}$ host given in \eqref{eqn:Dimless_Coupled_EclipDynB}, we obtain a system of equations, which is coupled to the ODE for the average density of virus in the air \eqref{eqn:AveBulk_Pathogen_3}, given by
\begin{align}\label{eqn:ReducedODE_EclipDyn}
\begin{split}
\frac{\text{d} T_j}{\text{d}t} &= -\beta_{1j}\,T_j  \left( V + \frac{\xi_j v_j}{2\pi D_0} \right) -  \frac{\mu}{D_0} \beta_{1j}\,T_j \,   \big{(} \mathcal{G} \Psi \big{)}_j   - \beta_{2j} T_j v_j, \\[1ex]
\frac{\text{d} E_{ j}}{\text{d}t} & =\beta_{1j}\,T_j  \left( V + \frac{\xi_j v_j}{2\pi D_0} \right) +  \frac{\mu}{D_0} \beta_{1j}\,T_j \,   \big{(} \mathcal{G} \Psi \big{)}_j   + \beta_{2j}  T_j v_j - k_j E_j,\\[1ex]
\frac{\text{d} I_{j}}{\text{d}t} & = k_j E_j - \delta_j I_j,\\[1ex]
\frac{\text{d} {v}_j}{\text{d} t} & = p_j {I}_j - c_j {v}_j, 
\end{split}
\end{align}
with initial conditions 
\begin{align}
\begin{split}
V(0) \ge 0, \quad T_{j}(0) &\ge 0, \quad E_{j}(0) \ge 0, \quad  I_{j}(0)\ge 0, \quad  v_{j}(0)\ge0, \quad \text{for} \quad j=1,\dots,m.
\end{split}
\end{align}  
\end{subequations}
This multiscale ODE system integrates the within-host dynamics of infected individuals to the spatial spread of virus between the hosts. It incorporates the effect of the diffusive viral particles and  host position on the within-host kinetics of the virus. The nonlinear ODE system is derived in the regime where the diffusion rate of  viral particles in the air is intermediate, as indicated in \eqref{eqn:rescale_D}.
 Comparing the ODE system to that of the classical TCL model \cite{ke2021vivo, perelson2021mechanistic, iyaniwura2024kinetics}, we observe that the difference between these two model come from the first two terms in the equation for target cells ($T_j$) and cells in eclipse phase ($E_j$). The first term in the $T_j$ and $E_j$ equations given by $\beta_{1j}\,T_j  \left( V + \frac{\xi_j v_j}{2\pi D_0} \right)$ describe the effect of the diffusing viral particles in the air on the within-host dynamics of the infected host, while the second term in these equations, given by $\frac{\mu}{D_0} \beta_{1j}\,T_j \,   \big{(} \mathcal{G} \Psi \big{)}_j$, models the effect of the location of the individuals in the domain on the within-host dynamics of the virus. The effect of the locations are captured through the Neumann Green's matrix $ \mathcal{G}$ defined in \eqref{eqn:G_matrix}.

\section{Existence and uniqueness of model solution}\label{sec:ExistUnique}

In this section, we analyze the mathematical properties of the multiscale ODE model~\eqref{eqn:ReducedODE_EclipDyn_IVP} to establish its well-posedness. Since the modeling framework presented here is novel, it is crucial to demonstrate that the system admits unique solutions for biologically relevant initial conditions. Proving the existence, uniqueness, and boundedness of solutions ensures that the model is both mathematically rigorous and biologically interpretable. 

To this end, we first establish the non-negativity and boundedness of the solutions, which are essential for biological realism. In particular, we present and prove Lemma~\ref{Lemma:NonNegative} to show that the model variables remain non-negative over time, reflecting the fact that cell populations and viral particle concentrations cannot be negative. Next, we prove Lemma~\ref{Lemma:Boundedness}, which guarantees that the solutions remain bounded for all finite times. Finally, we combine these results to demonstrate the existence and uniqueness of solutions to the ODE system~\eqref{eqn:ReducedODE_EclipDyn_IVP}, thereby confirming that the model is well-posed.

\begin{lemma}[Non-negativity]\label{Lemma:NonNegative}
    Given the initial conditions 
    \begin{equation}
    \begin{split}
    V(0) \ge 0, \quad T_{j}(0) &\ge 0, \quad E_{j}(0) \ge 0, \quad  I_{j}(0)\ge 0, \quad  v_{j}(0)\ge0, \quad \text{for} \quad j=1,\dots,m,
    \end{split}
    \end{equation}
    the solution $\Psi_j(t) = V(t)\cup \Phi_j(t) $, where $ \Phi_j(t) = \big{\{ } T_{j}(t), E_{j}(t), I_{j}(t), v_j(t) \big{\} }$ for $j=1,\dots,m,$
    of the multiscale model \eqref{eqn:ReducedODE_EclipDyn_IVP} is non-negative for all time $t > 0$.
\end{lemma}

\begin{proof}
We first show that $T_j(t) \geq 0$ for all $t>0$, for $j=1,\dots,m$. Suppose the initial number of target cells $T_j(0) \geq 0$. From the equation for $T_j(t)$ in \eqref{eqn:ReducedODE_EclipDyn}, we have
\begin{equation}\label{Eq:T_eqn_proof}
\frac{\text{d} T_{j}(t)}{dt}  = \bigg[- \beta_{1j} \,  \left( V(t) + \frac{\xi_j v_j(t)}{2\pi D_0} \right) -  \frac{\mu}{D_0} \beta_{1j} \,   \left(\mathcal{G}\, \Psi \right)_j   - \beta_{2j} v_j(t) \bigg] T_j.
\end{equation} 	
Applying the integrating factor method to \eqref{Eq:T_eqn_proof} and simplifying, we obtain 
\begin{align}\label{Eq:T_eqn_proof2}
T_{j}(t) & = T_{j}(0) \exp \bigg[- \int_{0}^{t} \beta_{1j} \,  \left( V(u) + \frac{\xi_j v_j(u)}{2\pi D_0} \right) du -  \frac{\mu}{D_0} \beta_{1j} \,   \int_{0}^{t} \Big( \mathcal{G}\Psi\Big)_j \, du   - \int_{0}^{t} \Big( \beta_{2j} v_j(u) \Big) du \bigg].
\end{align}
Since $T_j(0) \geq 0$, from \eqref{Eq:T_eqn_proof2}, we have  $T_{j}(t) \ge 0$ for all $t > 0$ and  for $j=1,\dots,m$.
Next, we show that $E_j(t) \geq 0$ for all $t>0$, for $j=1,\dots,m$. From the equation for $E_j$ in \eqref{eqn:ReducedODE_EclipDyn}, we have
\begin{equation}\label{Eq:E_eqn_proof}
\begin{split}    
\frac{\text{d} E_{j}(t)}{dt}  & = \beta_{1j}\,T_j  \left( V + \frac{\xi_j v_j}{2\pi D_0} \right) +  \frac{\mu}{D_0} \beta_{1j}\,T_j \,   \big{(} \mathcal{G} \Psi \big{)}_j   + \beta_{2j}  T_j v_j - k_j E_j .
\end{split}
\end{equation}
To continue with this proof, we need to show that 
\begin{equation}\label{Eq:T_greater_proof}
\bigg[ \beta_{1j} \left( V(t) + \frac{\xi_j v_j(t)}{2\pi D_0} \right) +  \frac{\mu}{D_0} \beta_{1j} \,   \left(\mathcal{G}\, \Psi \right)_j   + \beta_{2j} v_j(t) \bigg] \geq 0 \quad \text{for} \quad t>0.
\end{equation} 	
Since our model assumes that there is no replenishment of target cells during infection, the population of target cells ($T_j(t)$) in an infected individual will continue to decrease during infection. In other words, $T_j$ is monotone decreasing during infection, which implies that $\frac{\text{d} T_{j}(t)}{dt}  \leq 0$ for all $t >0$. Based on this, from \eqref{Eq:T_eqn_proof}, we have
\begin{equation}\label{Eq:T_greater_proof2}
\bigg[ \beta_{1j} \left( V(t) + \frac{\xi_j v_j(t)}{2\pi D_0} \right) +  \frac{\mu}{D_0} \beta_{1j} \,   \left(\mathcal{G}\, \Psi \right)_j   + \beta_{2j} v_j(t) \bigg] T_j(t) \geq 0.
\end{equation} 
Since $T_j(t) \geq 0$, we conclude that the inequality in \eqref{Eq:T_greater_proof} holds. Based on this, we therefore write \eqref{Eq:E_eqn_proof} as
\begin{equation}\label{Eq:E_eqn_proof11}
\begin{split}    
\frac{\text{d} E_{j}(t)}{dt}  & =  T_j(t) \left[ \beta_{1j} \left( V(t) + \frac{\xi_j v_j}{2\pi D_0} \right) +  \frac{\mu}{D_0} \beta_{1j} \,   \big{(} \mathcal{G} \Psi \big{)}_j   + \beta_{2j} \, v_j \right] - k_j E_j \ge - k_j E_j, 
\end{split}
\end{equation}
which implies that  
\begin{equation}\label{Eq:E_eqn_proof2}
 \frac{\text{d} E_{j}(t)}{dt} +  k_j E_j \geq 0  .
\end{equation}
Applying the integrating factor method to \eqref{Eq:E_eqn_proof2} and simplifying, we obtain
\begin{equation}
E_{j}(t) \geq E_{j}(0) \, e^{- k_j t}, \quad \text{for} \quad j=1,\dots,m,
\end{equation}
for all $t>0$ since $E_j(0) \geq 0$. 

To show that $I_j(t) \geq 0$ for all $t>0$, for $j=1,\dots,m$, we consider the equation for $I_j(t)$ in \eqref{eqn:ReducedODE_EclipDyn}. Since $E_j(t) \geq 0$ and $k_j \geq 0$, from this equation, we have
\begin{equation}\label{Eq:I_eqn_proof}
\frac{\text{d} I_{j}(t)}{dt} = k_j E_j - \delta_j I_j \ge - \delta_j I_j, 
\end{equation}
Integrating, we have
\begin{equation}
I_{j}(t) \geq I_{j}(0) \,e^{\delta_j \, t}  \geq 0,  \quad \text{for} \quad j=1,\dots,m,
\end{equation}
for all $t>0$ since the initial condition $I_j(0) \geq 0$ for $j=1,\dots,m$. We are left with showing the concentration of virus in the air $V(t)$ and the $j^{\text{th}}$ infected individual ($v_j$) are non-negative. 
From the equation for $v_j$ in \eqref{eqn:ReducedODE_EclipDyn}, we have
\begin{equation}\label{Eq:v_eqn_proof}
\frac{\text{d} {v}_j}{\text{d} t} = p_j {I}_j - c_j {v}_j(t)) \ge - c_j v_j(t),
\end{equation}
since $I_j(t) \geq 0$ and $p_j \geq 0$. Integrating and imposing the initial condition for $v_j$,  we obtain 
\begin{equation}
v_j (t) \ge v_j (0) \, e^{-c_j t} \geq 0,
\end{equation}
Therefore, $v_j(t) \geq 0$ for all $t>0$ and for $j=1,\dots,m$. Finally, we show the non-negativity of $V(t)$ for all $t>0$. 
From the equation for $V(t)$ in \eqref{eqn:ReducedODE_EclipDyn}, we have
\begin{equation}
 \frac{\text{d}V}{\text{d}t} + V = \frac{1}{|\Omega|} \sum_{i=1}^{m} \xi_i v_i \,.
\end{equation}
Integrating, gives
\begin{equation}\label{Eq:Vt_eqn_proof}
\begin{split}
V(t) & = \frac{e^{-t}}{|\Omega|} \,  \sum_{i=1}^{m} \left( \int_{0}^{t}  \xi_i \,v_i (u) \; \text{d}u \right)  + V(0)\, e^{-t}. 
\end{split}
\end{equation}
Since $v_i(t) \ge 0$ for $j=1,\dots,m$ as shown in \eqref{Eq:v_eqn_proof} and $V(0) \geq 0 $, we conclude that $V(t) \geq 0$ for all $t>0$. Thus, given non-negative initial conditions, the solutions of the multiscale model \eqref{eqn:ReducedODE_EclipDyn_IVP} is non-negative for all $t>0$.
\end{proof}

\begin{lemma}[Boundedness] \label{Lemma:Boundedness}
The closed set $\mathcal{D} = \cup_{j=1}^{m} \mathcal{D}_{j}$, where $\mathcal{D}_{j} = \Big.\Big\{\big(T_{j}, E_{j}, I_{j}\big) \in \mathbb{R}_{+}^{3}\; \Big| \;\mathcal{N}_{j}(t) \leq \mathcal{N}_{j}(0) \big.\big\}$ for $j=1,\dots,m$, and $\mathcal{N}_{j}(t) = T_{j} + E_{j} + I_{j}$, is positively invariant with respect to the multiscale model \eqref{eqn:ReducedODE_EclipDyn_IVP}. Furthermore, there exist constants $y \geq 0$ and $y_j \geq 0$ such that $V(t) \leq y$ and $v_j(t) \leq y_j$ for $t>0$ and $j=1,\dots,m$.    
\end{lemma}

\begin{proof}

Let $\mathcal{N}_{j}(t) = T_{j} + E_{j} + I_{j}$ for $j=1,\dots,m$. Differentiating $\mathcal{N}_{j}(t)$ with respect to $t$ and substituting the derivatives of $T_j(t), E_j(t)$, and $I_j(t)$ given in \eqref{eqn:ReducedODE_EclipDyn_IVP}, we obtain
\begin{equation}\label{Eq:Nt_deriv}
    \frac{\text{d} \mathcal{N}_{j}(t)}{dt}  = - \delta_{j} I_{j} \leq 0 \, , 
\end{equation}
since $I_j(t) \geq 0$ and $\delta_j \geq 0$. This implies that $\mathcal{N}_{j}(t)$ is monotone decreasing for all $t > 0$. Combining this with the fact that $T_j(t), E_j(t)$, and $I_j(t)$  are non-negative  for all $t > 0$, as shown in Lemma~\ref{Lemma:NonNegative}, it suffices that  $\mathcal{N}_{j}(t) \to 0$ as $t \to \infty$. We need to show that $\mathcal{N}_{j}(t)$ is bounded above. Integrating \eqref{Eq:Nt_deriv}, we get
\begin{equation}\label{Eq:N_Integral}
~  \int_{0}^{t} d \mathcal{N}_{j}(s) \leq K,
\end{equation}
where $K$ is some constant. Since $\mathcal{N}_{j}(t)$ is monotone decreasing, any non-negative value of  $K$ satisfies \eqref{Eq:N_Integral} for $t >0 $. Without loss of generality, we set $K=0$. So that
\begin{equation}
 \mathcal{N}_{j}(t) \leq \mathcal{N}_{j}(0),    \quad \text{for} \quad j=1,\dots,m.
\end{equation}
This implies that $\mathcal{N}_{j}(t)$ is bounded above by the initial values of $\mathcal{N}_{j}$. Since $\mathcal{N}_{j}(t) = T_{j}(t) + E_{j}(t) + I_{j}(t)$ is monotone decreasing, and $\mathcal{N}_{j}(t) \to 0$ as $t \to \infty$, and the quantities $T_j(t), E_j(t)$, and $I_j(t)$ are non-negative, we conclude that $T_j(t), E_j(t)$, and $I_j(t)$ are bounded for all $t > 0$ and $j=1,\dots,m$.
Thus, $\mathcal{D} = \cup_{j=1}^{m} \mathcal{D}_{j}$  is positively invariant with respect to the multiscale model \eqref{eqn:ReducedODE_EclipDyn_IVP}.

Next, we show that $v_j$ is bounded. We have already showed that $v_j(t) \geq 0 $ for all $t>0$ in Lemma~\ref{Lemma:NonNegative}, which shows that $v_j$ is bounded below. We are left to show that it is bounded above.
From the equation for $v_j$ in \eqref{eqn:ReducedODE_EclipDyn_IVP}, we have 
\begin{equation}\label{Eq:v_eqn_proof_2}
\frac{\text{d} {v}_j}{\text{d} t} = p_j {I}_j(t) - c_j  \,v_j \leq  p_j {I}_j(t),
\end{equation}
since $v_j \geq 0$ and $I_j \geq 0$,  and the parameters $p_j$, and $c_j$ are all positive. Integrating, we have
\begin{equation}
\begin{split}
v_j (t) & \leq v_j(0) + p_j \, \int_{0}^{t} I_j(u) \; \text{d}u.
\end{split}
\end{equation}
We know that the initial number of infected cells is zero, i.e., $I_j(0) = 0$ and $I_j(t) \to 0$ as $t \to \infty$, which implies that $I(t)$ is a concave function. We assume that $I_j(t)$ approaches zero rapidly as $t \to \infty$ such that the integral $\int_{0}^{t} I_j(u) \; \text{d}u$ is bounded. In this regard, there exist a constant $y_j \geq 0$ such that 
\begin{equation}\label{Eq:vj_Bnded}
\begin{split}
v_j (t) & \leq v_j(0) + p_j \, \int_{0}^{t} I_j(u) \; \text{d}u \leq y_j\, .
\end{split}
\end{equation}
This implies that $v_j (t)  \leq y_j$ for all $t> 0$ and $j=1,\dots,m$.

Lastly, we need to show that $V(t)$ is bounded. Consider the equation for $V(t)$ in \eqref{eqn:ReducedODE_EclipDyn_IVP}, from which we obtain
\begin{equation}
 \frac{\text{d}V(t)}{\text{d}t} = \frac{1}{|\Omega|} \sum_{i=1}^{m} \xi_i \,v_i(t) - V(t) \leq \frac{1}{|\Omega|} \sum_{i=1}^{m} \xi_i \, v_i(t),
\end{equation}
since $\xi_j \geq 0$ and $v_j(t) \geq 0$ and $V(t) \geq 0$  for all $t >0$. Integrating, we have 
\begin{equation}
 V(t) \leq V(0) + \frac{1}{|\Omega|} \sum_{i=1}^{m} \xi_i \left( \int_{0}^{t}  \,v_i(u) \; \text{d}u \right).
\end{equation}
Similar to $I_j(t)$, the dynamics of virus $(v_j(t))$ in infected individuals is concave. In addition, since $v_j \geq 0$, we assume that $v_j$ approaches zero rapidly as $t \to \infty$, so that the integral 
$\int_{0}^{t}  \,v_j(u) \; \text{d}u$ is bounded. Therefore, there exists a constant $y \geq 0$ such that
\begin{equation}\label{Eq:V_Bnded}
 V(t) \leq V(0) + \frac{1}{|\Omega|} \sum_{i=1}^{m} \xi_i \left( \int_{0}^{t}  \,v_i(u) \; \text{d}u \right) \leq y,
\end{equation}
which implies that $ V(t) \leq y$ for all $y \geq 0$.
\end{proof}

Now, we have all the tools we need to prove the existence of a unique solution to the multiscale model \eqref{eqn:ReducedODE_EclipDyn_IVP}.
We state the following result based on the classical Picard–Lindelöf theorem, which ensures the existence and uniqueness of solutions to initial value problems under suitable conditions on continuity and Lipschitz continuity \cite{walter2013ordinary, coddington1955theory, hale2009ordinary, siegmund2016generalized, teschl2012ordinary}.

\begin{theorem}[Existence and uniqueness]\label{Thm:ExistUnique}
The multiscale ODE model 
\begin{subequations}\label{eqn:IVP_theorem}
    \begin{equation}\label{eqn:ODE_theorem}
    \begin{split}
         \frac{\text{d}\,V(t)}{\text{d}t} &= h_1(V, v_j);\\[1ex]
         \frac{\text{d}\,Y_{kj}(t)}{\text{d}t} &= h_{kj}(V,  T_j, E_j, I_j, v_j); \quad \text{for} \quad k=2,\dots,5 \quad \text{and}  \quad j=1,\dots,m,
    \end{split}
\end{equation}
where $(Y_{2j}, Y_{3j}, Y_{4j}, Y_{5j}) \equiv ( T_j, E_j, I_j, v_j)$ and the functions $h_1$ and $h_{kj}$ are given by
\begin{align}\label{eqn:h_Functions}
\begin{split}
h_1 &= \frac{1}{|\Omega|} \sum_{i=1}^{m} \xi_i v_i - V, \\[1ex]
h_{2j} &= -\beta_{1j}\,T_j  \left( V + \frac{\xi_j v_j}{2\pi D_0} \right) -  \frac{\mu}{D_0} \beta_{1j}\,T_j \,   \big{(} \mathcal{G} \Psi \big{)}_j   - \beta_{2j} T_j v_j, \\[1ex]
h_{3j}  & =\beta_{1j}\,T_j  \left( V + \frac{\xi_j v_j}{2\pi D_0} \right) +  \frac{\mu}{D_0} \beta_{1j}\,T_j \,   \big{(} \mathcal{G} \Psi \big{)}_j   + \beta_{2j}  T_j v_j - k_j E_j,\\[1ex]
h_{4j}  & = k_j E_j - \delta_j I_j,\\[1ex]
h_{5j} & = p_j {I}_j - c_j {v}_j, 
\end{split}
\end{align}
\end{subequations}
with  initial conditions 
\begin{align}
\begin{split}
 V(0) \geq 0, \quad 
  T_j(0) & \geq 0, \quad  E_j(0) \geq 0, \quad  I_j(0) \geq 0, \quad v_j(0) \geq 0, \qquad  j=2,\dots,m, 
\end{split}
\end{align}   
has a unique solution, which is continuously differentiable for $t \in [0, \infty)$ if the following conditions are satisfied:
\begin{enumerate}
    \item The functions $h_1$ and $h_{kj}$ for $k=2,\dots,5$ and $j=1,\dots,m$ in Eq.~\eqref{eqn:h_Functions}, defined on a domain \(D = [0, \infty) \times \mathbb{R}^5\) for each $j$, are continuous in both  $t$ and the variables $V,  T_j, E_j, I_j$, and $v_j$ for $j=1,\dots,m$.
    \item The partial derivatives of the functions $h_1$ and $h_{kj}$ for $k=2,\dots,5$ and $j=1,\dots,m$, with respect to the variables $V,  T_j, E_j, I_j$, and $v_j$ are bounded on the domain $D$ for $j=1,\dots,m$.
\end{enumerate}
\end{theorem}

\begin{proof}
 Since the functions $h$, $h_{1j}$, $h_{2j}$, $h_{3j}$, and $h_{4j}$, as given in \eqref{eqn:h_Functions} are all polynomial functions in the variables $V$, $T_j$, $E_j$, $I_j$, and $v_j$, it can easily be verified that these functions are continuous in these variables and $t$. 

For the second part of this proof, we need to show that the following partial derivatives are bounded:
\begin{equation}\label{Eq:JacMatrix}
\frac{\partial h_1}{\partial Z_k}, \quad   \frac{\partial h_{2j}}{\partial Z_k}, \quad   \frac{\partial h_{3j}}{\partial Z_k}, \quad   \frac{\partial h_{4j}}{\partial Z_k},
\quad  \text{and} \quad \frac{\partial h_{5j}}{\partial Z_k}, \quad \textbf{for} \quad k = 1, \dots, 5 \quad \text{and} \quad j=1,\dots, m,
\end{equation}
where $Z_k$ is the $k^{\text{th}}$ entry of the vector $Z = (V, T_j, E_j, I_j, v_j)$.
The majority of these partial derivatives are zero. 

For $h_1(V, T_j, E_j, I_j, v_j)$, we have $ \frac{\partial h_{1}}{\partial V} = -1$ and  $\frac{\partial h_{1}}{\partial T_j} = \frac{\partial h_{1}}{\partial E_j}= \frac{\partial h_{1}}{\partial I_j} = 0$ for $j=1,\dots, m$. Lastly, we have
\begin{equation}
    \frac{\partial h_{1}}{\partial v_j} = \frac{\xi_j}{|\Omega|}.
\end{equation}
Therefore, all the partial derivatives of $h_1(V, T_j, E_j, I_j, v_j)$ with respect to the model variables are bounded. Similarly, for $h_{2j}(V, T_j, E_j, I_j, v_j)$, we have $ \frac{\partial h_{2j}}{\partial E_j} =  \frac{\partial h_{2j}}{\partial I_j} = 0$ and
\begin{equation}
    \left| \frac{\partial h_{2j}}{\partial V} \right|  = \left| -\beta_{1j}\,T_j \right| \le  \beta_{1j}\, \mathcal{N}_{j}(0). 
\end{equation}
Differentiating $h_{2j}$ with respect to $T_j$ and using the fact that \eqref{Eq:T_greater_proof} holds, and triangular inequality, we have
\begin{equation}
    \begin{split}
  \left| \frac{\partial h_{2j}}{\partial T_{j}} \right| & <  \beta_{1j} \left| V \right| +  \frac{\beta_{1j}\,\xi_j}{2\pi D_0} \left| v_j  \right| +  \frac{\mu \beta_{1j} \xi_1 }{D_0} \,   \left|R_{1} \right| \left| v_1 \right| +  \frac{\mu \beta_{1j} \xi_j }{D_0} \sum_{i=2}^{m}  \left|G_{1,i}\right| \left|v_i\right|  +  \beta_{2j} \left|v_j\right|.
\end{split}
\end{equation}
Using the boundedness properties of $v_j$ and $V$ in \eqref{Eq:vj_Bnded} and  \eqref{Eq:V_Bnded}, respectively, we conclude that $\left| \frac{\partial h_{2j}}{\partial T_{j}} \right| $ is bounded. Lastly, we differentiate  $h_{2j}$ with respect to $v_j$ and use triangular inequalities on the resulting derivative to obtain
\begin{equation}
    \begin{split}
  \left| \frac{\partial h_{2j}}{\partial v_{j}} \right| & \leq  \frac{\beta_{1j}\,\xi_j}{2\pi D_0} \left|T_j \right|  + \frac{\mu \beta_{1j} \xi_1}{D_{0}} \left| T_{j}\right| \left|R_{1}\right| +  \frac{\mu \beta_{1j} \xi_j }{D_0}  \left| T_{j}\right| \sum_{i=2}^{m} \xi_i   \left|G_{1,i}\right|   + \beta_{2j} \left|T_j \right|, \\
  & \leq  \left[ \frac{\beta_{1j}\,\xi_j}{2\pi D_0}   + \frac{\mu \beta_{1j} \xi_1}{D_{0}} \left|R_{1}\right|  + \sum_{i=2}^{m} \xi_i   \left|G_{1,i}\right|   + \beta_{2j} \right] \mathcal{N}_{j}(0).
\end{split}
\end{equation}
Therefore, we conclude that the partial derivatives of $h_{2j}$ with respect to all the model variables are bounded.

Next, we differentiate $h_{3j}(V, T_j, E_j, I_j, v_j)$ with respect to $V$ and $E_j$ to obtain
\begin{equation}
    \left| \frac{\partial h_{3j}}{\partial V} \right| =  \beta_{1j} \left| T_j \right| \le  \beta_{1j}\, N_j (0) \qquad \text{and} \qquad  \left| \frac{\partial h_{3j}}{\partial E_{j}} \right|  =  k_j \;,
\end{equation}
respectively. Note that $\left| \frac{\partial h_{3j}}{\partial I_j} \right| = 0$. For $T_j$, we have
\begin{equation}\label{Eq:h3_proof1}
\begin{split}
\left|\frac{\partial h_{3j}}{\partial T_{j}}\right| & = \left| \beta_{1j} \left( V + \frac{\xi_j v_j}{2\pi D_0} \right) +  \frac{\mu}{D_0} \beta_{1j}\,  \left( \mathcal{G} \, \psi \right)_{j}   + \beta_{2j} v_j \right|.
\end{split}    
\end{equation}
Since the inequality in \eqref{Eq:T_greater_proof} holds, using triangular inequality on \eqref{Eq:h3_proof1}, we obtain
\begin{equation}
\begin{split}
\left|\frac{\partial h_{3j}}{\partial T_{j}}\right| & < \beta_{1j} \left| V \right| +  \frac{\beta_{1j}\,\xi_j}{2\pi D_0} \left| v_j  \right| +  \frac{\mu \beta_{1j} \xi_1 }{D_0} \,   \left|R_{1} \right| \left| v_1 \right| +  \frac{\mu \beta_{1j} \xi_j }{D_0} \sum_{i=2}^{m}  \left|G_{1,i}\right| \left|v_i\right|  +  \beta_{2j} \left| v_j\right|,
\end{split}    
\end{equation}
which implies that $\left| \frac{\partial h_{3j}}{\partial T_{j}} \right| $ is bounded since  $v_j$ and $V$ are bounded for $t > 0$. Differentiating $h_{3j}$ with respect to $v_j$, we have
\begin{equation}
\begin{split}
\left| \frac{\partial h_{3j}}{\partial v_{j}} \right| <  \left[ \frac{\beta_{1j}\,\xi_j}{2\pi D_0} +  \frac{\mu \beta_{1j} \xi_1 }{D_0} \,   \left|R_{1} \right|  +  \frac{\mu \beta_{1j} \xi_j }{D_0} \sum_{i=2}^{m}  \left|G_{1,i}\right|   +  \beta_{2j} \right] N_j(0).
\end{split}
\end{equation}
Therefore, the partial derivatives of $h_{3j}$ with respect to all the model variables are bounded. Lastly, for the functions $h_{4j}$ and $h_{5j}$, we have 
\begin{equation}
\begin{split}
 \left| \frac{\partial h_{4j}}{\partial E_{j}} \right| = k_j, \quad \left| \frac{\partial h_{4j}}{\partial I_{j}} \right| = \delta_j, \quad
\left| \frac{\partial h_{5j}}{\partial I_{j}} \right|  = p_j, \quad \text{and} \quad \left|\frac{\partial h_{5j}}{\partial v_{j}}\right| =  c_j ,\\ 
\end{split}
\end{equation}
and the remaining derivatives are zeros. This shows that the derivatives of  $h_{4j}$ and $h_{5j}$ for $j=1,\dots,m$ are all bounded. \end{proof}

\section{Analysis of within-host viral kinetics}\label{sec:ViralProperties}

Now that we have established the existence, uniqueness, and boundedness of the solutions of our multiscale ODE model \eqref{eqn:ReducedODE_EclipDyn_IVP}, we want to use this model to study the within-host dynamics of the disease. First, we consider the ODE model \eqref{eqn:ReducedODE_EclipDyn_IVP} in different limits. At leading-order, we have the ODE system which is independent of the host locations, given by:
\begin{equation}\label{eqn:ReducedODE_EclipDyn_IVP_Lead}
\begin{split}
    \frac{\text{d}V(t)}{\text{d}t} &= \frac{1}{|\Omega|} \sum_{i=1}^{m} \xi_i\, v_i(t) - V(t), \\
\frac{\text{d} T_j}{\text{d}t} &= -\beta_{1j}\,T_j  \left( V + \frac{\xi_j v_j}{2\pi D_0} \right)   - \beta_{2j} T_j v_j, \\[1ex]
\frac{\text{d} E_{ j}}{\text{d}t} & =\beta_{1j}\,T_j  \left( V + \frac{\xi_j v_j}{2\pi D_0} \right)   + \beta_{2j}  T_j v_j - k_j E_j,\\[1ex]
\frac{\text{d} I_{j}}{\text{d}t} & = k_j E_j - \delta_j I_j,\\[1ex]
\frac{\text{d} {v}_j}{\text{d} t} & = p_j {I}_j - c_j {v}_j.
\end{split}
\end{equation}

We observe from \eqref{eqn:ReducedODE_EclipDyn_IVP_Lead}, that in the limit of fast diffusing viral particles in the air, i.e.,  $D_0 \to \infty$, we have $\frac{\xi_j v_j}{2\pi D_0} \to 0$. Therefore, the ODE system reduces to
\begin{equation}\label{eqn:ReducedODE_EclipDyn_IVP_WM}
\begin{split}
    \frac{\text{d}V(t)}{\text{d}t} &= \frac{1}{|\Omega|} \sum_{i=1}^{m} \xi_i\, v_i(t) - V(t), \\
\frac{\text{d} T_j}{\text{d}t} &= -\beta_{1j}\,T_j  V    - \beta_{2j} T_j v_j, \\[1ex]
\frac{\text{d} E_{ j}}{\text{d}t} & =\beta_{1j}\,T_j V    + \beta_{2j}  T_j v_j - k_j E_j,\\[1ex]
\frac{\text{d} I_{j}}{\text{d}t} & = k_j E_j - \delta_j I_j,\\[1ex]
\frac{\text{d} {v}_j}{\text{d} t} & = p_j {I}_j - c_j {v}_j,
\end{split}
\end{equation}
which represents the well-mixed version of the ODE model \eqref{eqn:ReducedODE_EclipDyn_IVP}. By ignoring the contribution of viral particles in the air (i.e., setting \( V = 0 \)), the well-mixed model reduces to the classical target cell limited (TCL) model \cite{perelson2002modelling, perelson2021mechanistic, ke2021vivo, iyaniwura2024kinetics}, which describes the in vivo infection kinetics for a host independent of the viral pathogens expelled by the remaining hosts, and is given by:
\begin{equation}\label{eqn:TCL_Model}
\begin{split}
\frac{\text{d} T_j}{\text{d}t} &=    - \beta_{2j} T_j v_j, \\[1ex]
\frac{\text{d} E_{ j}}{\text{d}t} & =   \beta_{2j}  T_j v_j - k_j E_j,\\[1ex]
\frac{\text{d} I_{j}}{\text{d}t} & = k_j E_j - \delta_j I_j,\\[1ex]
\frac{\text{d} {v}_j}{\text{d} t} & = p_j {I}_j - c_j {v}_j.
\end{split}
\end{equation}   
This derivation demonstrates that the classical TCL model emerges as a limiting case of our multiscale framework when the airborne diffusion of viral particles is neglected. In this case, the infection kinetics of the hosts are decoupled.

Next, we study the within-host properties of the ODE model \eqref{eqn:ReducedODE_EclipDyn_IVP} and compare them with those of the classical TCL model \eqref{eqn:TCL_Model}. We consider a scenario with only one host in the domain, and another where there are multiple hosts. 
Since the model has only one steady state (the virus-free steady state), we focus our analysis on the basic reproduction number of the model and the sensitivity of some infection kinetics metric to the model parameters.

\subsection{Model with a single host}\label{subsec:R0_Single}

For a single host, 
the  ODE model \eqref{eqn:ReducedODE_EclipDyn_IVP} reduces to 
\begin{align}\label{eqn:SingleInd_Model}
\begin{split}
\frac{\text{d}V}{\text{d}t} &= \frac{\xi\, v(t)}{|\Omega|}  - V(t), \\[1ex]
\frac{\text{d} T}{\text{d}t} &= -\beta_{1}\,T  \left( V + \frac{\xi v}{2\pi D_0} \right) -  \frac{\mu}{D_0} \beta_{1}\,T \,   \big{(} \mathcal{G} \Psi \big{)}   - \beta_{2} T v, \\[1ex]
\frac{\text{d} E}{\text{d}t} & =\beta_{1}\,T  \left( V + \frac{\xi v}{2\pi D_0} \right) +  \frac{\mu}{D_0} \beta_{1}\,T \,   \big{(} \mathcal{G} \Psi \big{)}   + \beta_{2}  T v - k E,\\[1ex]
\frac{\text{d} I}{\text{d}t} & = k E - \delta I,\\[1ex]
\frac{\text{d} {v}}{\text{d} t} & = p {I} - c \,{v}, 
\end{split}
\end{align}
where $\big{(} \mathcal{G} \Psi \big{)} = R \, \xi v$ with $R \equiv R(\pmb{x})$ defined as the singular part of the Neumann Green's function, given by \eqref{Eq:GreensFunction_Sing}.
We compute the within-host basic reproduction number of the virus for an individual, which gives the number of cells a single viral particle will infect in a newly infected individual. 
We use the next-generation matrix approach \cite{van2002reproduction}, where the next-generation matrix is given by $M = \mathbb{F}\mathbb{V}^{-1}$, with $\mathbb{F}$ defined as the matrix of new infections and $\mathbb{V}$ as the matrix for the transfer of infection.

In terms of the model in Eq.~\eqref{eqn:SingleInd_Model}, at the virus-free equilibrium, VFE $\equiv \left(V_e, T_e, E_e, I_e, v_e \right) = \left(0, \mathcal{N}(0), 0, 0, 0 \right)$, where $\mathcal{N}(0) = T(0) + E(0) + I(0) + v(0)$, the matrices $\mathbb{F}$ and $\mathbb{V}$ are given by
\begin{displaymath}
	\mathbb{F} = \begin{pmatrix}
	0 & 0 & 0 & 0 \\
	\\
	\beta_{1} \mathcal{N}(0)& 0 & 0 & \beta_{1}  \mathcal{N}(0) \xi  \left(\frac{1 + 2 \pi  \mu  R}{2 \pi  D_0} \right) + \beta_{2} \mathcal{N}(0) \\
	\\
	0 & 0 & 0 & 0 \\
	\\
	0 & 0 & 0 & 0 \\
	\end{pmatrix} \quad \text{and} \quad
    \mathbb{V} = \begin{pmatrix}
	 1 & 0 & 0 & -\frac{\xi}{|\Omega|} \\
	 \\
	 0 & k & 0 & 0 \\
	 \\
	 0 & - k & \delta & 0 \\
	 \\
	 0 & 0 & -p & c\\ 
	\end{pmatrix} .
\end{displaymath}
Using these matrices to compute the next-generation matrix ($M$), we obtain our desire basic reproduction number, which is the spectral radius of the matrix.
The two-term asymptotic expansion of the virus within-host basic reproduction number in a newly infected host, $\mathcal{R}_{0}$, is given by
\begin{equation}\label{Eq:R0_SingVIralParticle}
\mathcal{R}_{0} = \mathcal{R}_{1}  + \frac{ \mu }{D_0} \mathcal{R}_{2} \,+  \mathcal{O}\Bigg( \left( \frac{ \mu }{D_0} \right)^2 \Bigg),
\end{equation}
where $\mathcal{R}_{1}$ and $\mathcal{R}_{2}$ are defined by
\begin{equation}\label{Eq:R0_SingInd_NoLoc}   
\mathcal{R}_{1} = \frac{p\,\beta_{2}  \,   \mathcal{N}(0)}{\delta c}  +  \frac{p  \,\beta_{1}   \xi  \,  \mathcal{N}(0) \Big(2 \pi D_0 + |\Omega| \Big)}{ 2 \pi \delta D_0 |\Omega|   c  } \quad
\text{and} \quad \mathcal{R}_{2} = \frac{  p  \, \beta_{1} \, \xi R \, \mathcal{N}(0)}{ \delta  c }.
\end{equation}
As in the ODE model \eqref{eqn:SingleInd_Model}, the leading-order term in the basic reproduction number ($\mathcal{R}_{1})$ incorporates the effect of the domain geometry and diffusion rate of particles in the air into $\mathcal{R}_0$, while the $\mathcal{O}(1)$ correction term incorporates the effect of the host location.
As $D_0 \to \infty$, the basic reproduction number in Eq.~\eqref{Eq:R0_SingVIralParticle} reduces to
\begin{equation}\label{Eq:R0_SingVIralParticle_WellMixed}
\mathcal{R}_{0} = \frac{p\,\beta_{2}  \,   \mathcal{N}(0)}{\delta c } + \frac{ p \, \beta_{1}  \,  \xi   \mathcal{N}(0) }{\delta  c  |\Omega|},
\end{equation}
which is the corresponding basic reproduction number for the well-mixed ODE model in Eq.~\eqref{eqn:ReducedODE_EclipDyn_IVP_WM}. Observe that this quantity is influenced by the size of the domain through the second term. As the domain becomes very large, i.e., $|\Omega| \to \infty$, the second term in \eqref{Eq:R0_SingVIralParticle_WellMixed} vanishes, and basic reproduction number reduces to
\begin{equation}\label{Eq:R0_SingInd_VDModel}
\mathcal{R}_{0} = \frac{p\,\beta_{2}  \,   \mathcal{N}(0)}{\delta c}.
\end{equation}
This implies that the domain has no effect of the within-host dynamics of the host. The corresponding basic reproduction number is that of the classical TCL model shown in Eq.~\eqref{eqn:TCL_Model}. As we have shown earlier, the classical TCL model is a limiting case of our multiscale model in the regime of fast diffusion of viral particles in the environment, and when each host is isolated with no spatial effect. Here, we just showed that the within-host basic reproduction number for the TCL model is a limiting case of the within-host reproduction number for our multiscale model in the well-mixed limit as $D_0 \to \infty$.

\subsubsection{Sensitivity analysis}\label{sec:Sensitiviy_analysis}

In this section, we assess the sensitivity of the model outputs to variations in important model parameters. We employed a global sensitivity analysis approach using Latin Hypercube Sampling (LHS) combined with Partial Rank Correlation Coefficients (PRCC) \cite{Blower,mckay2000comparison,stein1987large,loh1996latin}.

We began by specifying biologically reasonable lower and upper bounds for the model parameters. Using Latin Hypercube Sampling (LHS), we generated a large number of parameter sets to ensure uniform and stratified coverage of the parameter space. Each parameter set was used to run the model numerically 1000 times in MATLAB version R2025a \cite{MATLAB}. For each simulation, we recorded two key outcomes: the within-host basic reproduction number and the peak viral load. To account for potential non-linear but monotonic relationships, both the model outputs and input parameters were rank-transformed prior to analysis.

\begin{figure}[h]
  \centering
  \includegraphics[scale=0.22]{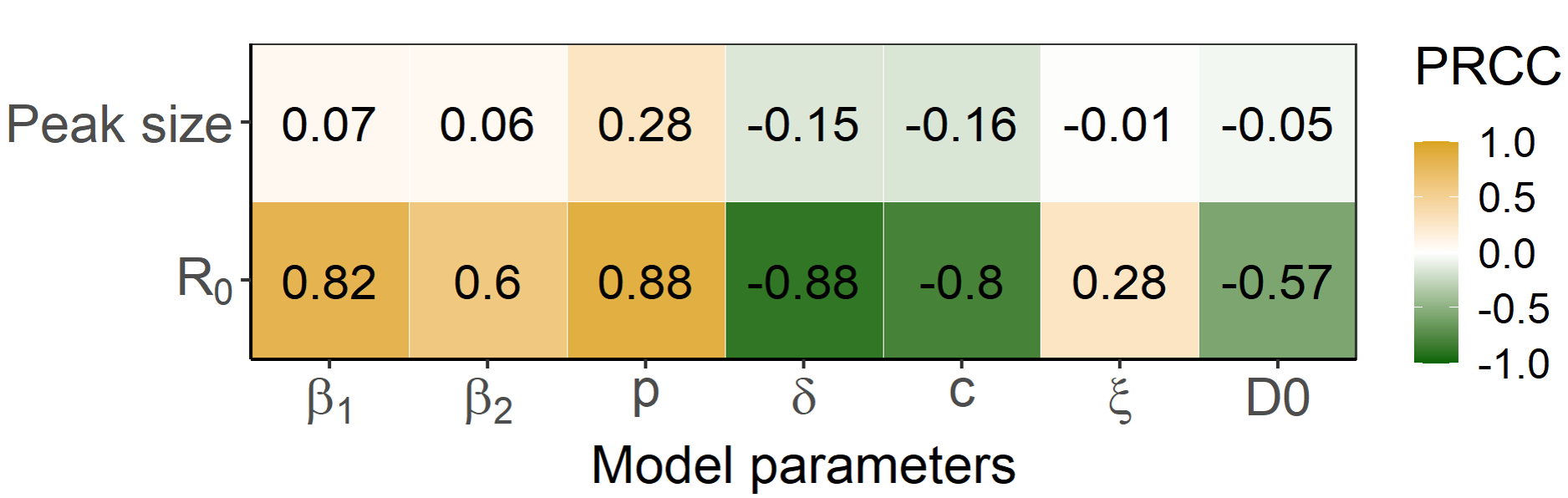}
  \caption{\textbf{Sensitivity analysis for the multiscale ODE model.} Partial Rank Correlation Coefficients (PRCC) for the basic reproduction number $(\mathcal{R}_0)$ in Eq.~\eqref{Eq:R0_SingVIralParticle} and viral load peak size for the multiscale model in Eq.~\eqref{eqn:SingleInd_Model} with respect to selected model parameters.}
\label{Fig:PRCC_DiffModel}
\end{figure}

Next, we conducted linear regression analyses to obtain the residuals for each model parameter and the associated output response. The Partial Rank Correlation Coefficients (PRCCs) were then calculated as the Pearson correlations between the residuals of the ranked inputs and outputs, providing a measure of both the strength and direction of each parameter’s effect on the outcome. These PRCC values were then visualized to identify the parameters with the greatest influence on infection kinetics (Figure~\ref{Fig:PRCC_DiffModel}).

\begin{figure}[h]
  \centering
  \includegraphics[scale=0.20]{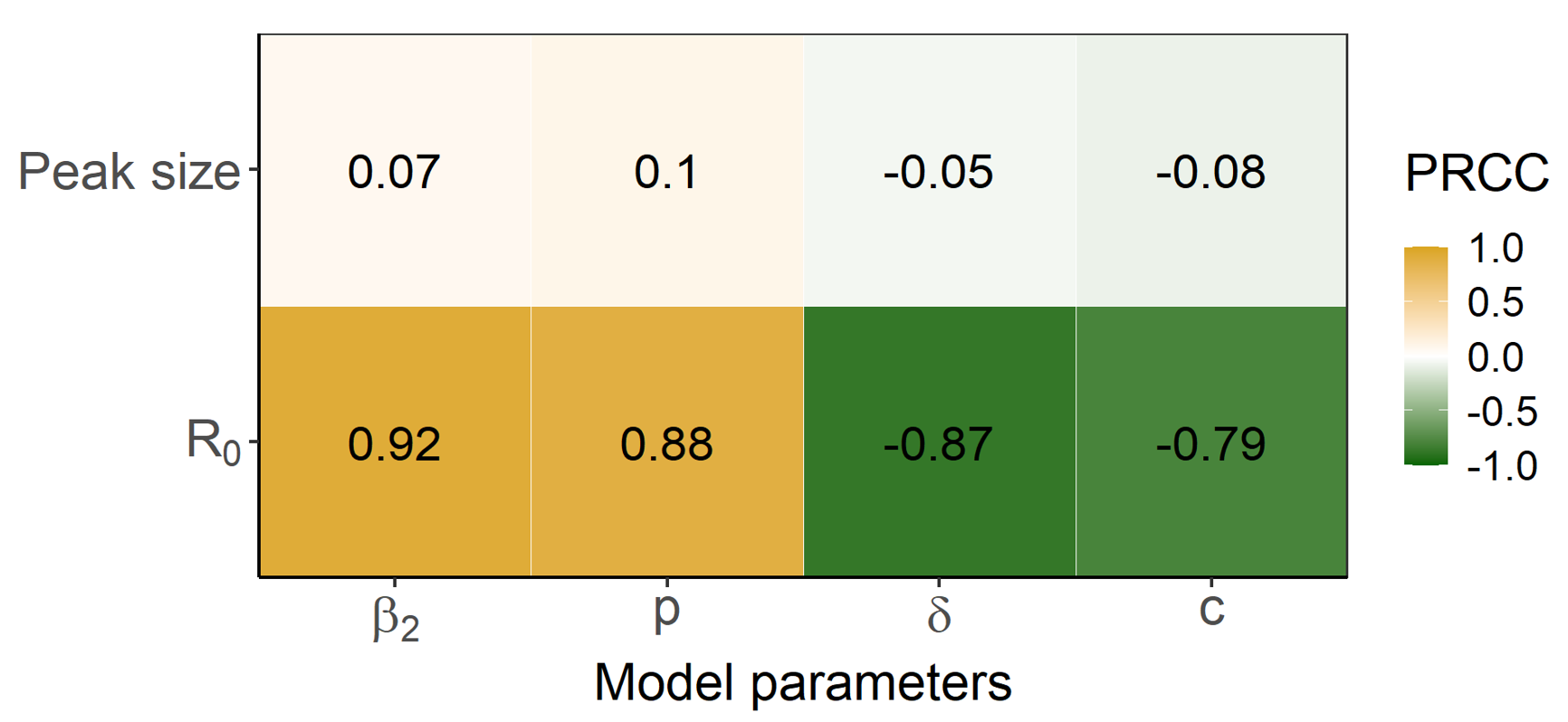}
  \caption{\textbf{Sensitivity analysis for the TCL model.} Partial Rank Correlation Coefficients (PRCC) for the basic reproduction number $(\mathcal{R}_0)$ in Eq.~\eqref{Eq:R0_SingInd_VDModel} and viral load peak size for the TCL model in Eq.~\eqref{eqn:TCL_Model} with respect to selected model parameters.}
\label{Fig:PRCC_TCLModel}
\end{figure}

We conduct a similar sensitivity analysis for the TCL model~\eqref{eqn:TCL_Model}, with results shown in Figure~\ref{Fig:PRCC_TCLModel}. In both models, the infection transmission rate ($\beta_2$) and viral production rate ($p$) are positively correlated with the basic reproduction number and viral peak size, indicating that increases in these parameters amplify disease spread and viral burden. Conversely, the death rate of infected cells ($\delta$) and the viral clearance rate ($c$) are negatively correlated with these outcomes, suggesting that higher values mitigate viral replication and peak viral load.

For the multiscale model~\eqref{eqn:SingleInd_Model}, additional parameters influence disease dynamics. The transmission rate due to diffusing virus ($\beta_1$) and the virus expulsion rate by infected hosts ($\xi$) are both positively correlated with $\mathcal{R}_0$, while their effects on viral peak size differ: $\beta_1$ increases the peak, whereas $\xi$ reduces it. Finally, the viral diffusion rate in air ($D_0$) is negatively correlated with both $\mathcal{R}_0$ and viral peak size, indicating that increased diffusion reduces disease transmission and viral accumulation.

\begin{figure}[h]
  \centering
  \includegraphics[width=1.65in]{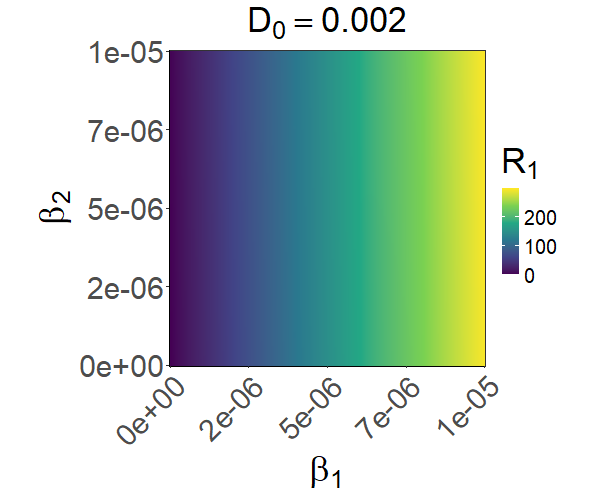}
  \includegraphics[width=1.65in]{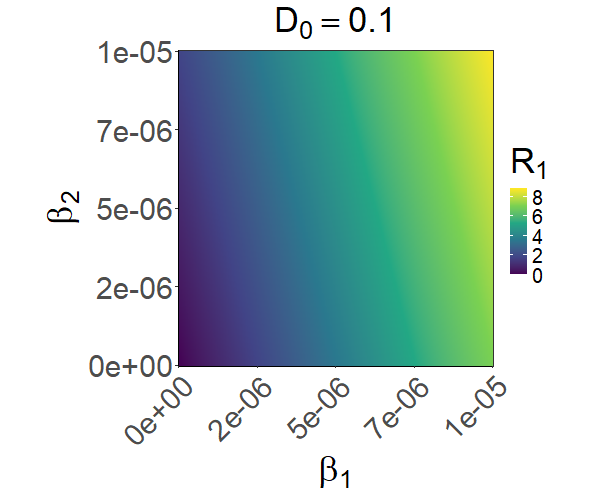}
  \includegraphics[width=1.65in]{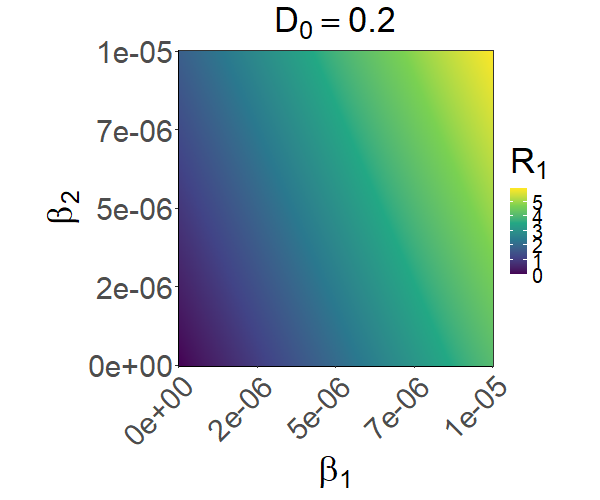}
   \includegraphics[width=1.65in]{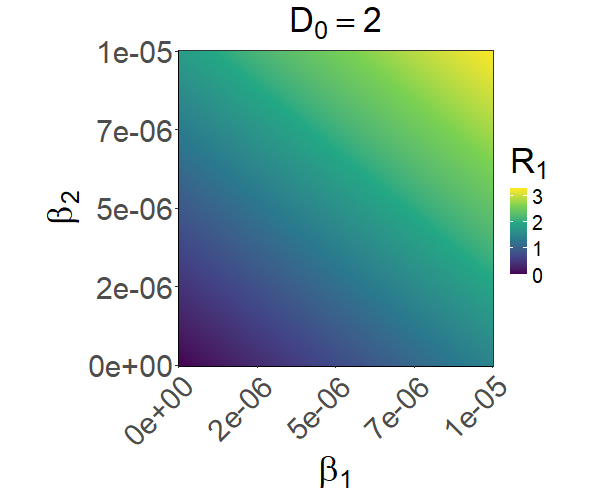}\\[2ex]
  \includegraphics[width=1.65in]{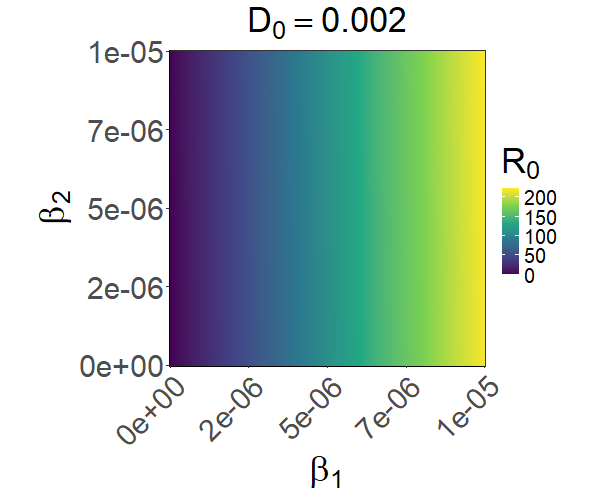}
  \includegraphics[width=1.65in]{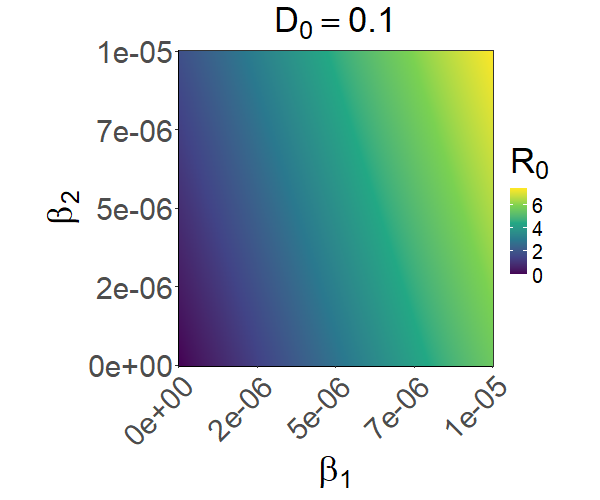}
  \includegraphics[width=1.65in]{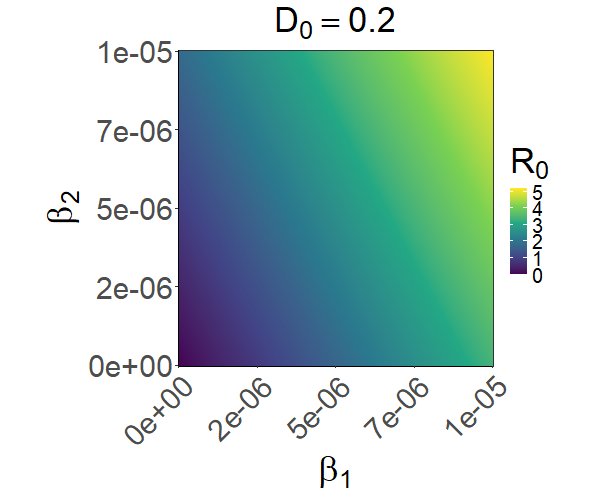}
   \includegraphics[width=1.65in]{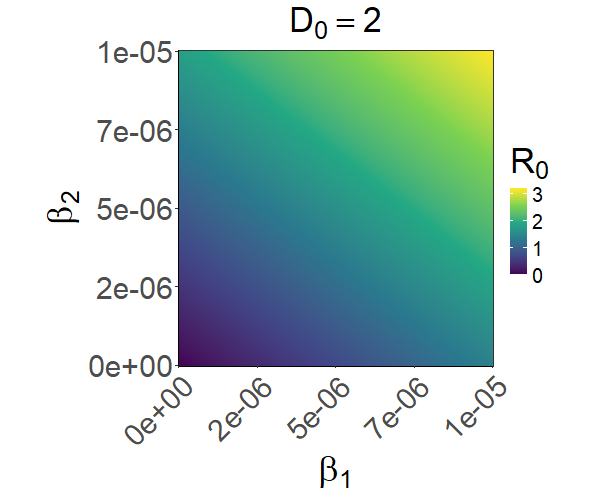}
\caption{\textbf{Effect of diffusion on within-host infection and the basic reproduction number.} 
Contour plots of the leading-order approximation of the basic reproduction number ($\mathcal{R}_1$) given in Eq.~\eqref{Eq:R0_SingInd_NoLoc} (top row) and the two-term expansion of the basic reproduction number ($\mathcal{R}_0$) in Eq.~\eqref{Eq:R0_SingVIralParticle} (bottom row). Results are shown as functions of the transmission rates $\beta_1$ and $\beta_2$ for different diffusion rates, ranging from $D_0 = 0.002$ to $D_0 = 2$. Parameters are given in Table~\ref{Table:TableParameter}.}
\label{Fig:geometry}
\end{figure}

In the next analysis, we investigate the impact of the viral diffusion rate on the basic reproduction number. To this end, we generate contour plots of the leading-order approximation of the basic reproduction number ($\mathcal{R}_1$), given in Eq.~\eqref{Eq:R0_SingInd_NoLoc}, and the two-term expansion of the reproduction number ($\mathcal{R}_0$) in Eq.~\eqref{Eq:R0_SingVIralParticle}, as functions of the transmission rates $\beta_1$ and $\beta_2$ for diffusion rates ranging from $D_0 = 0.002$ to $D_0 = 2$ (Figure~\ref{Fig:geometry}). The results in the top and bottom rows of Figure~\ref{Fig:geometry} differ only slightly, indicating that the leading-order term dominates $\mathcal{R}_0$, while the $\mathcal{O}(\mu/D_0)$ correction has a relatively small contribution. Moreover, for all values of $D_0$ considered, the maximum value of $\mathcal{R}_1$ exceeds that of $\mathcal{R}_0$, despite the absence of an explicit negative term in the expression for $\mathcal{R}_0$. This reduction arises because the singular component of the Neumann Green’s function ($R$) in Eq.~\eqref{Eq:GreensFunction_Sing} is negative in this setting, thereby lowering the overall value of $\mathcal{R}_0$.

\begin{table}[!h]
\caption{\textbf{Model parameters, descriptions, and values.}  }
\centering
\begin{tabular}{c  p{10.0cm}  >{\raggedleft\arraybackslash}p{3.250cm} >{\raggedleft\arraybackslash}p{2.0cm}}
\hline
{\bf{Par}} & {\textbf{Description}} &  {\textbf{Value}}  &  {\textbf{References}}  \\ \hline
$b_{1}\, (\beta_{1})$ &  Dimensional (dimensionless) infectious rate of diffusing virus & $10^{-8}$ ($5.6 \times 10^{-7}$) & \cite{wang2020modeling, ke2021vivo} \\[0.5ex]
$b_{2}\, (\beta_{2})$ &  Dimensional (dimensionless) infectious rate of in vivo virus & $10^{-8}$ ($5.6 \times 10^{-7}$) & \cite{wang2020modeling, ke2021vivo} \\[0.5ex]
$\gamma \, (\xi)$ & Dimensional (dimensionless) virus expulsion  rate & 2 (4.19) & Assumed \\[0.5ex]
$\alpha \, (k)$ & Dimensional (dimensionless) transition rate from eclipse cells to productively infected cells & 4 (0.4) & \cite{ke2021vivo, iyaniwura2024kinetics}  \\
$d \, (\delta)$ & Dimensional (dimensionless) death rate of infected cells  &  1.7 (0.17) & \cite{pawelek2012modeling, ke2021vivo} \\[0.5ex]
$\rho \, (p)$ & Dimensional (dimensionless) production rate of viral particles by infected cells  &  $10^{2}$ ($1.6 \times 10^{10}$)  & \cite{ke2021vivo, iyaniwura2024kinetics}   \\[0.5ex]
$D_r \, (D_0)$ & Dimensional (dimensionless) diffusion rate of viral particles in the air & Varied & Assumed \\[0.5ex]
$\phi \, (c)$ & Dimensional (dimensionless) within-host virus clearance rate  & 10 (1) & \cite{ke2021vivo, iyaniwura2024kinetics} \\[0.5ex]
$k_{r}$ & Dimensional degradation rate of viral particles in the air  & 10  & Assumed \\[0.5ex]
$r_c$ &  Typical value for the density of virus in the air   & 1 & Assumed \\[0.5ex]
$\varepsilon$ & Dimensionless radius of the disk used to represent individuals & 0.05 &  Assumed\\[0.5ex]
$|\Omega|$ & Area of the domain representing the enclosed environment in which the individuals are located & $\pi$  & Assumed  \\[0.5ex]
$N(0)$ & Total number of target cells & $8 \times 10^{7}$  & \cite{ke2021vivo, iyaniwura2024kinetics}  \\ 
 \hline
\end{tabular}
\label{Table:TableParameter}
\end{table}

For small diffusion rates, both $\mathcal{R}_1$ and $\mathcal{R}_0$ are nearly insensitive to $\beta_2$ and increase primarily with $\beta_1$, indicating that infection dynamics are dominated by diffusing viral particles entering the host. As $D_0$ increases, the relative contributions of $\beta_1$ and $\beta_2$ shift, with the influence of $\beta_2$ becoming more pronounced. At higher diffusion rates, both transmission pathways contribute comparably to the basic reproduction number, as illustrated by the case $D_0 = 2$. 
Moreover, the reproduction number is larger for small values of $D_0$ and decreases as $D_0$ increases. This trend is consistent with the results of the sensitivity analysis, which show that $\mathcal{R}_0$ is negatively correlated with $D_0$ (see Figure~\ref{Fig:PRCC_DiffModel}). Finally, the contour plots for $\mathcal{R}_1$ and $\mathcal{R}_2$ are nearly identical when $D_0 = 2$, suggesting that the system approaches a well-mixed regime at this diffusion rate.

\begin{figure}[!h]
  \centering
  \includegraphics[width=2.0in]{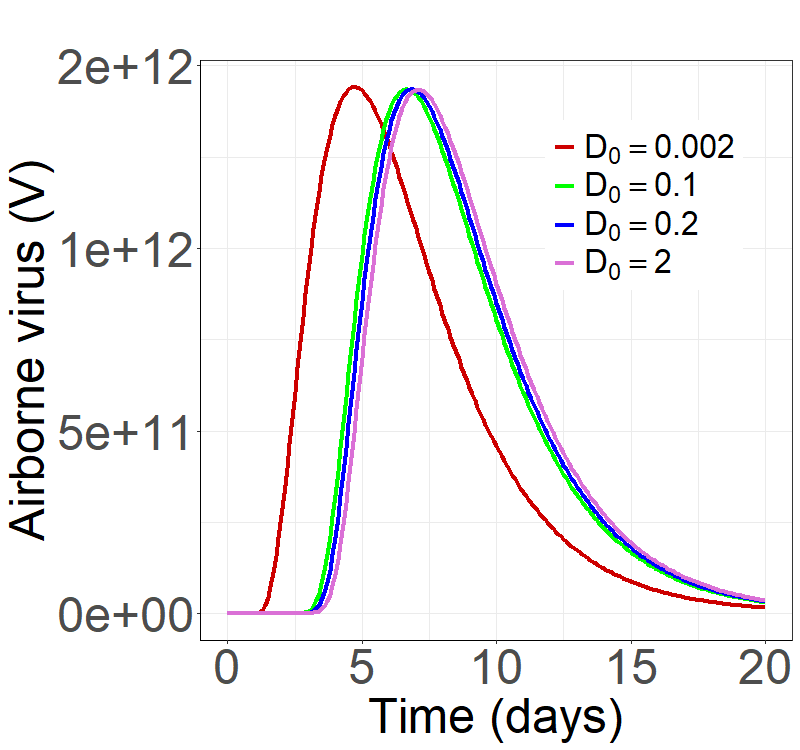} \quad
  \includegraphics[width=2.0in]{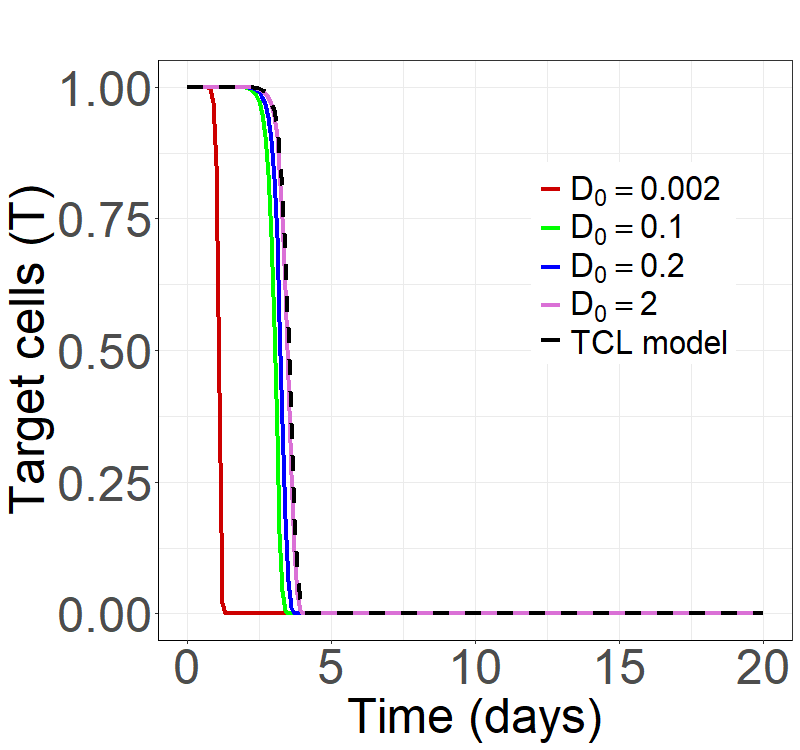} \quad
  \includegraphics[width=2.0in]{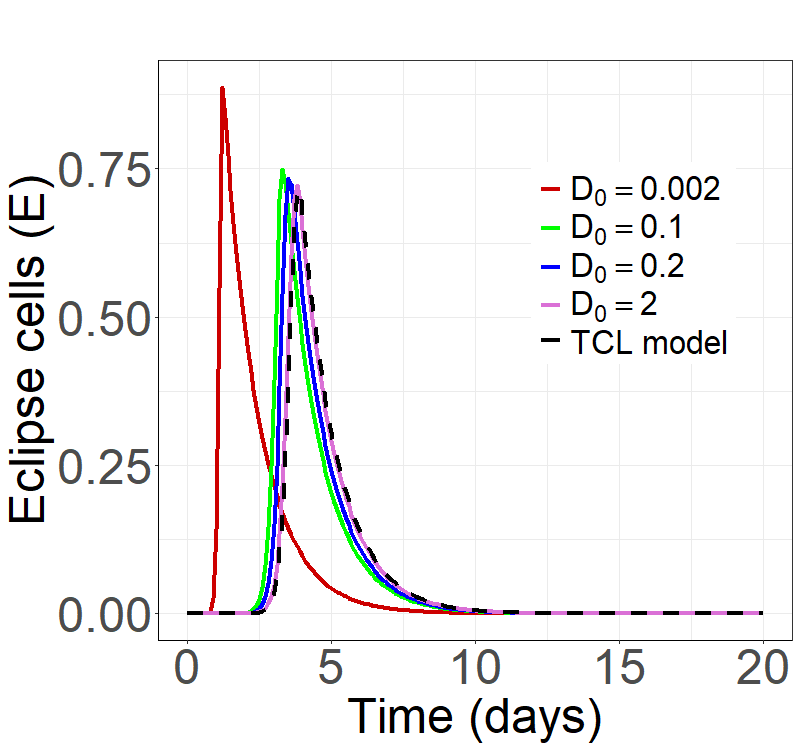}\\
  \includegraphics[width=2.0in]{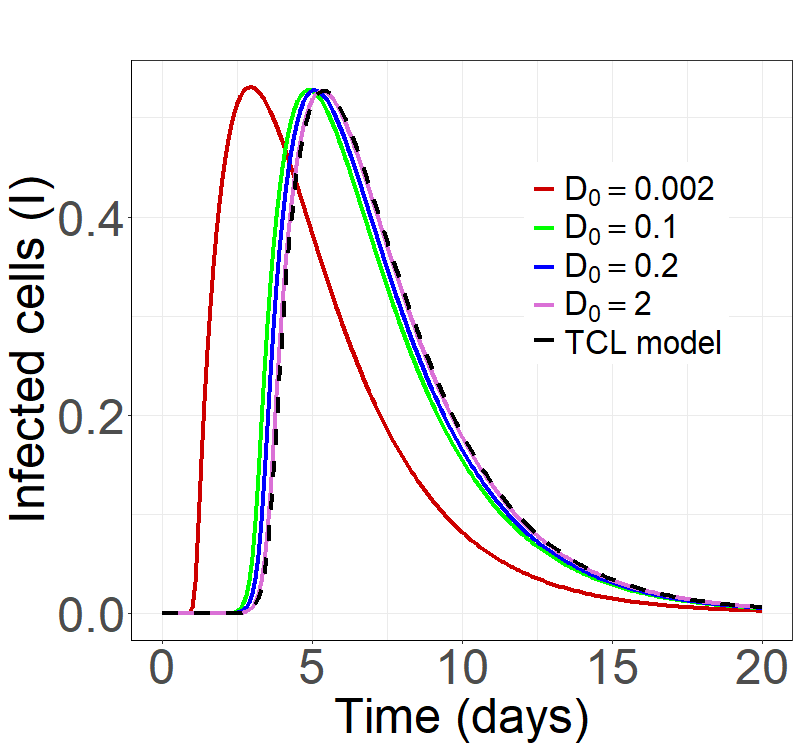} \quad
  \includegraphics[width=2.0in]{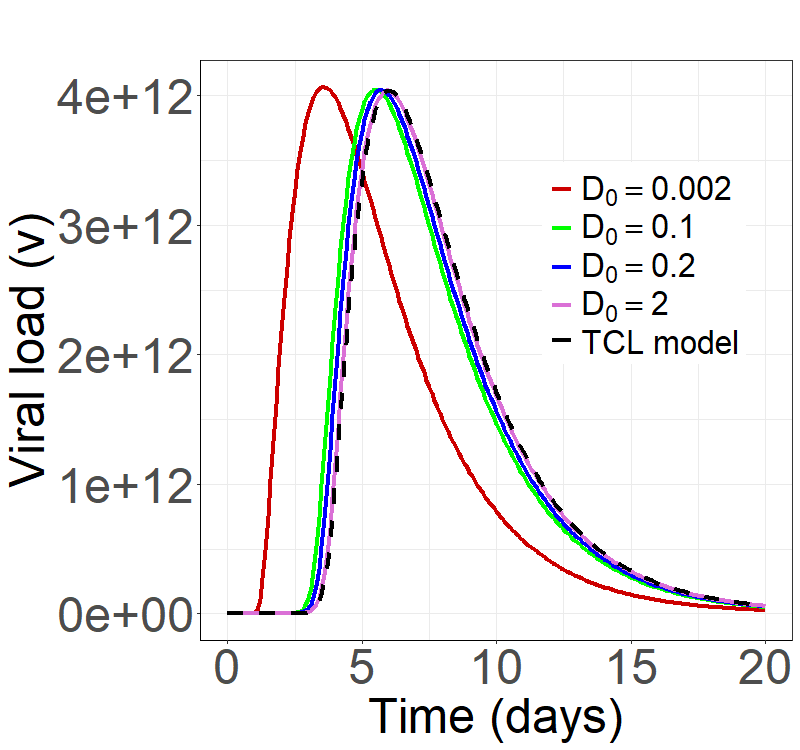}
  \caption{\textbf{Effect of diffusion on within-host infection dynamics.}
Time series of airborne viral concentration and within-host infection kinetics for a single host located at the origin $(0,0)$, computed using the multiscale model \eqref{eqn:SingleInd_Model} (solid lines) and the TCL model \eqref{eqn:TCL_Model} (dashed lines). The figure illustrates how diffusion influences within-host infection dynamics. Model parameters are provided in Table~\ref{Table:TableParameter}.}
\label{Fig:SingleHostDynamics}
\end{figure}

Next, we investigate the effect of diffusion on within-host infection dynamics in a simplified setting consisting of a single individual in the domain. Specifically, we examine how variations in the diffusion rate of airborne viral particles influence the dynamics of airborne virus, target cells, eclipse-phase cells, infected cells, and within-host viral particles, as computed by the multiscale model \eqref{eqn:SingleInd_Model}. These dynamics are compared with those obtained from the classical TCL model \eqref{eqn:TCL_Model} (see Figure~\ref{Fig:SingleHostDynamics}).
In these simulations, infection is initiated by a small number of cells in the eclipse phase. In particular, we set $E(0) = 1/N(0)$ for both models, with initial conditions $T(0)=1$, $V(0)=0$, $I(0)=0$, and $v(0)=0$. Parameter values used in both models are listed in Table~\ref{Table:TableParameter}.

From Figure~\ref{Fig:SingleHostDynamics}, we observe an initial delay before infection develops within the host. Such delays have been reported in previous studies and are often attributed to stochastic effects during the early stages of infection \cite{iyaniwura2024kinetics}. As the diffusion rate increases, this delay becomes more pronounced in the multiscale model. This behavior arises because higher diffusion disperses viral particles more rapidly throughout the domain, reducing their local concentration near the infected host and consequently decreasing infection via inhaled airborne virus.
As $D_0$ continues to increase, the solutions of the multiscale model converge to those of the TCL model. In particular, for $D_0 = 2$, we observe strong agreement between the two models. This finding is consistent with the contour plots shown in the last column of Figure~\ref{Fig:geometry}, which indicate that the system is effectively well-mixed at this diffusion rate.

\subsection{Model with multiple hosts}

\begin{figure}[!h]
  \centering
  \includegraphics[width=2.0in]{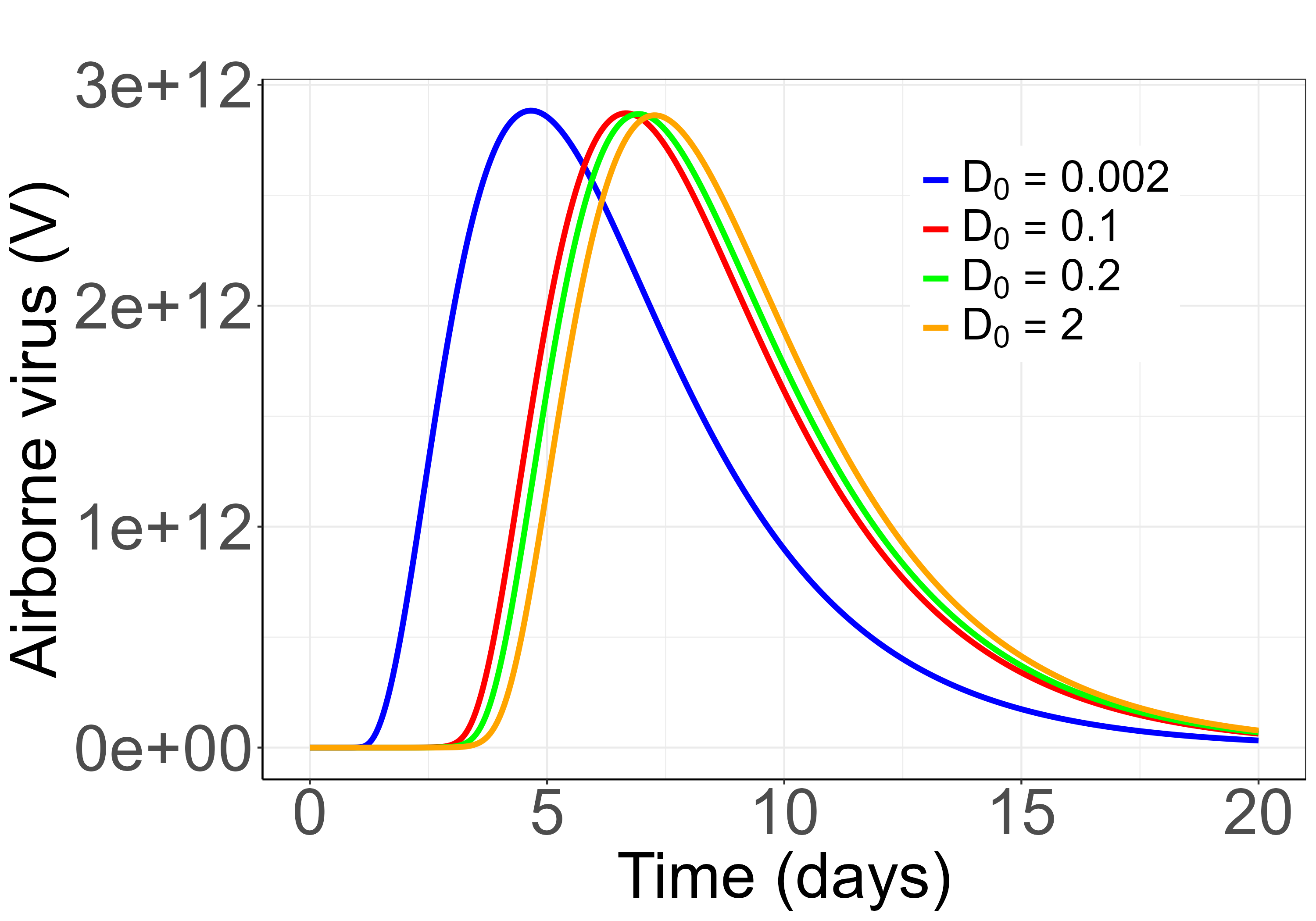} \quad
  \includegraphics[width=2.0in]{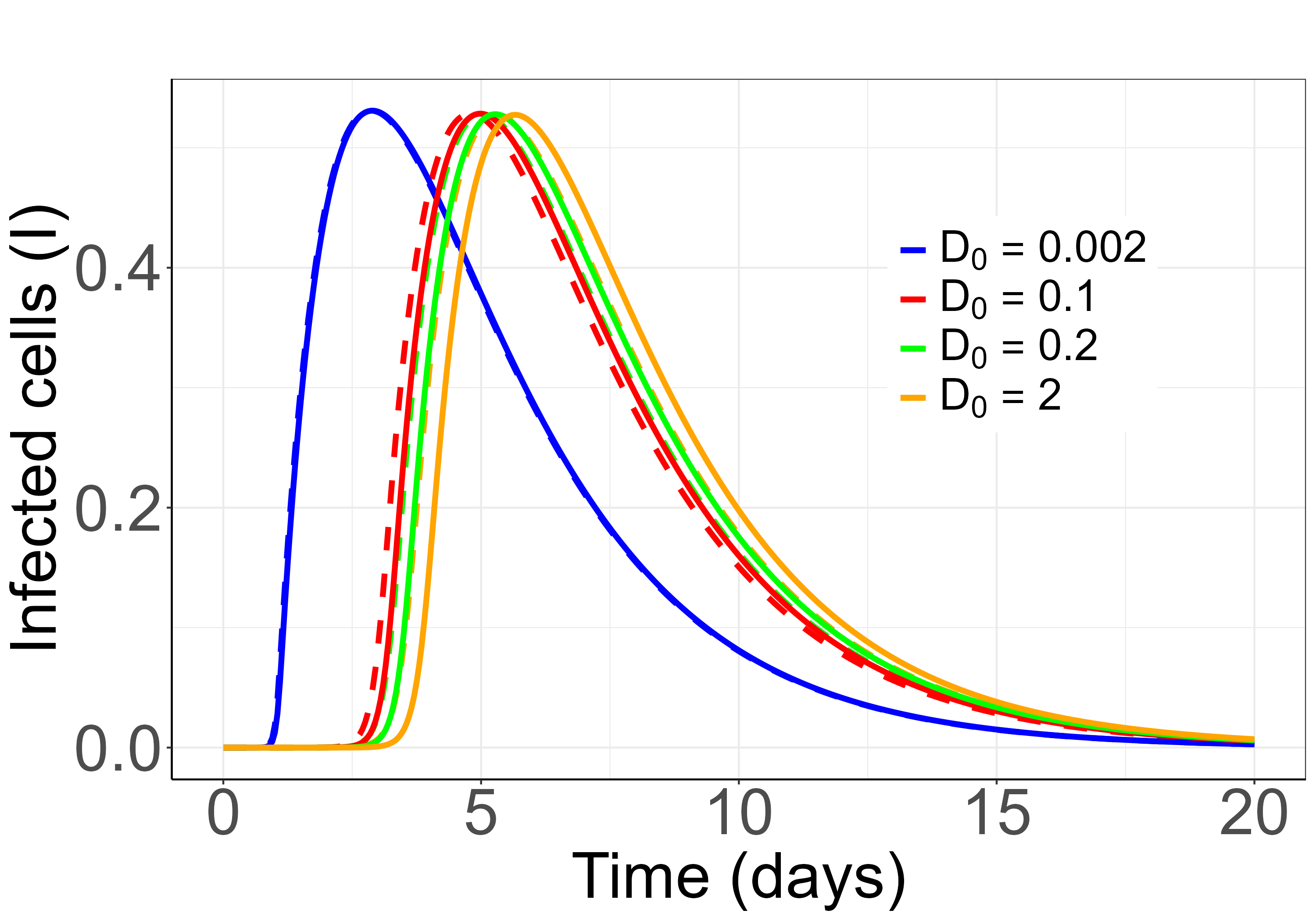} \quad
  \includegraphics[width=2.0in]{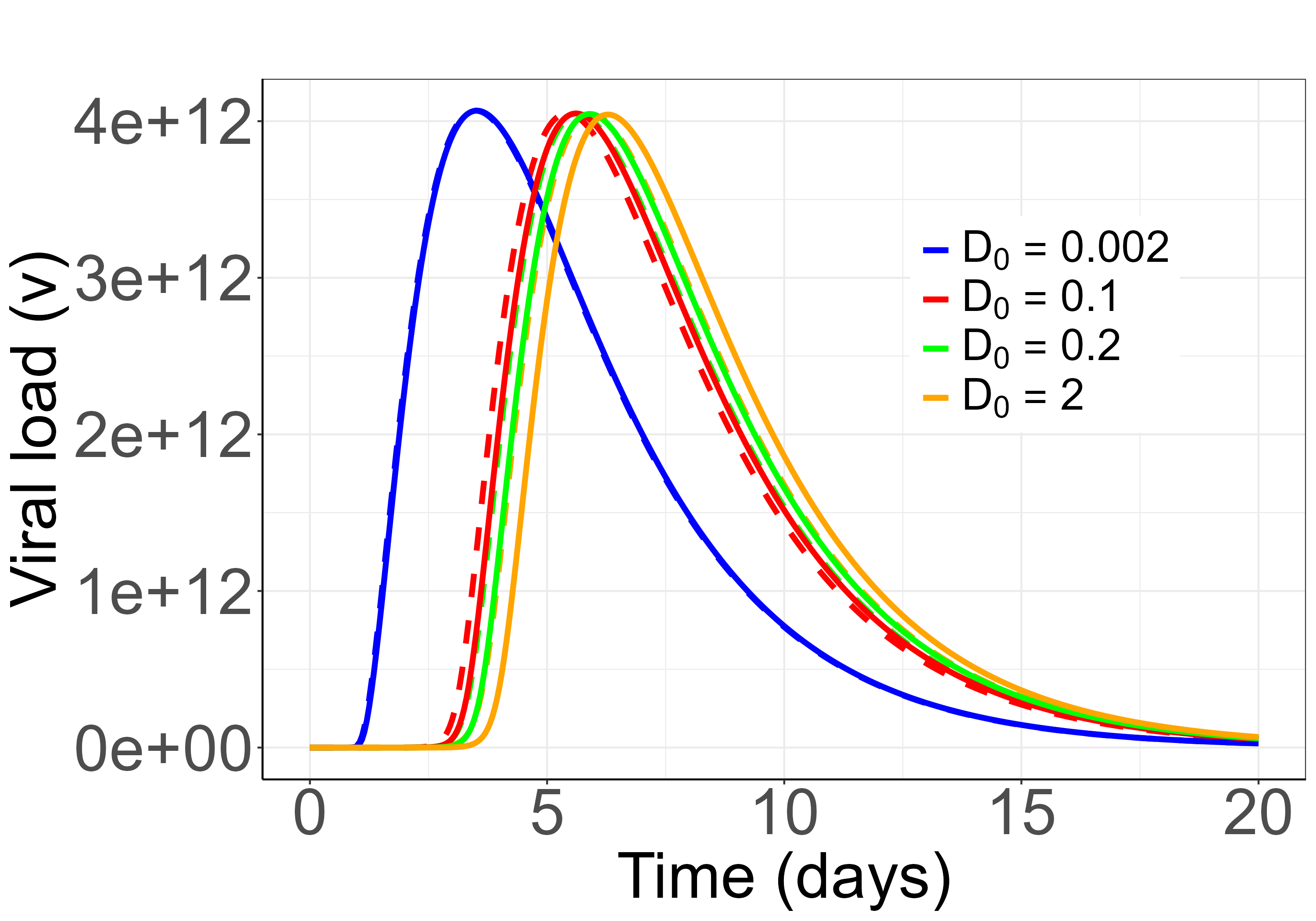}\\
  \includegraphics[width=2.0in]{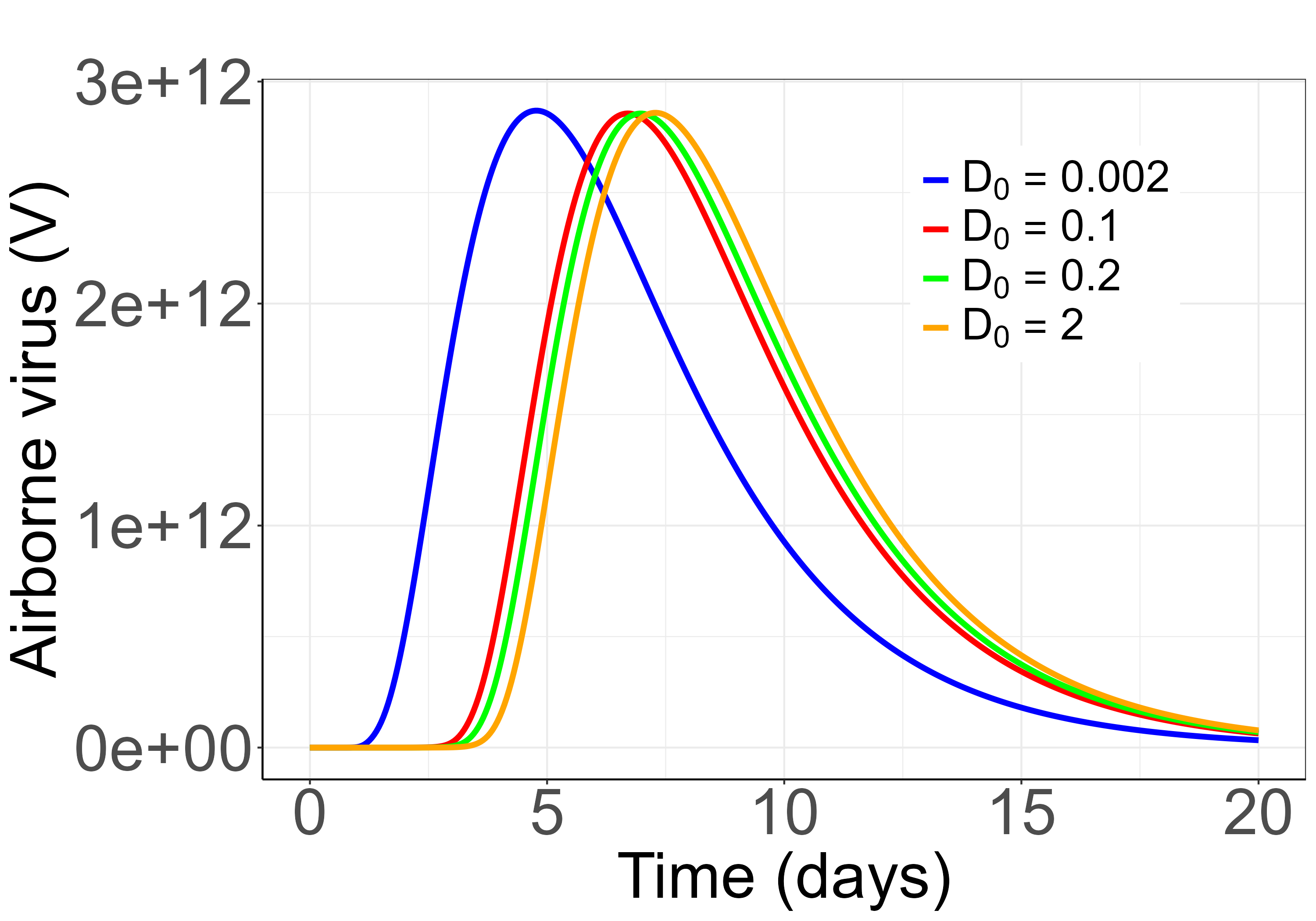} \quad
  \includegraphics[width=2.0in]{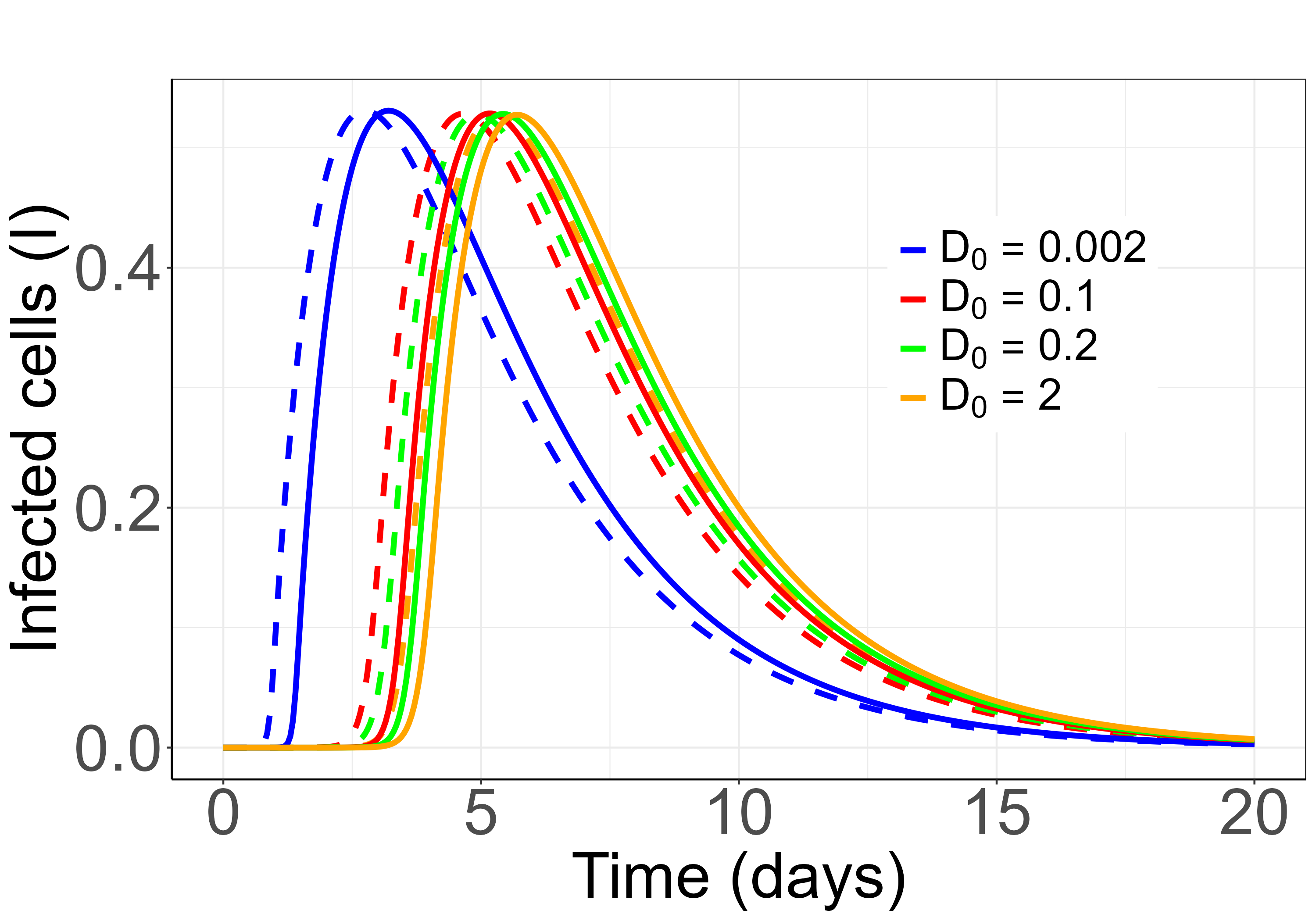} \quad
  \includegraphics[width=2.0in]{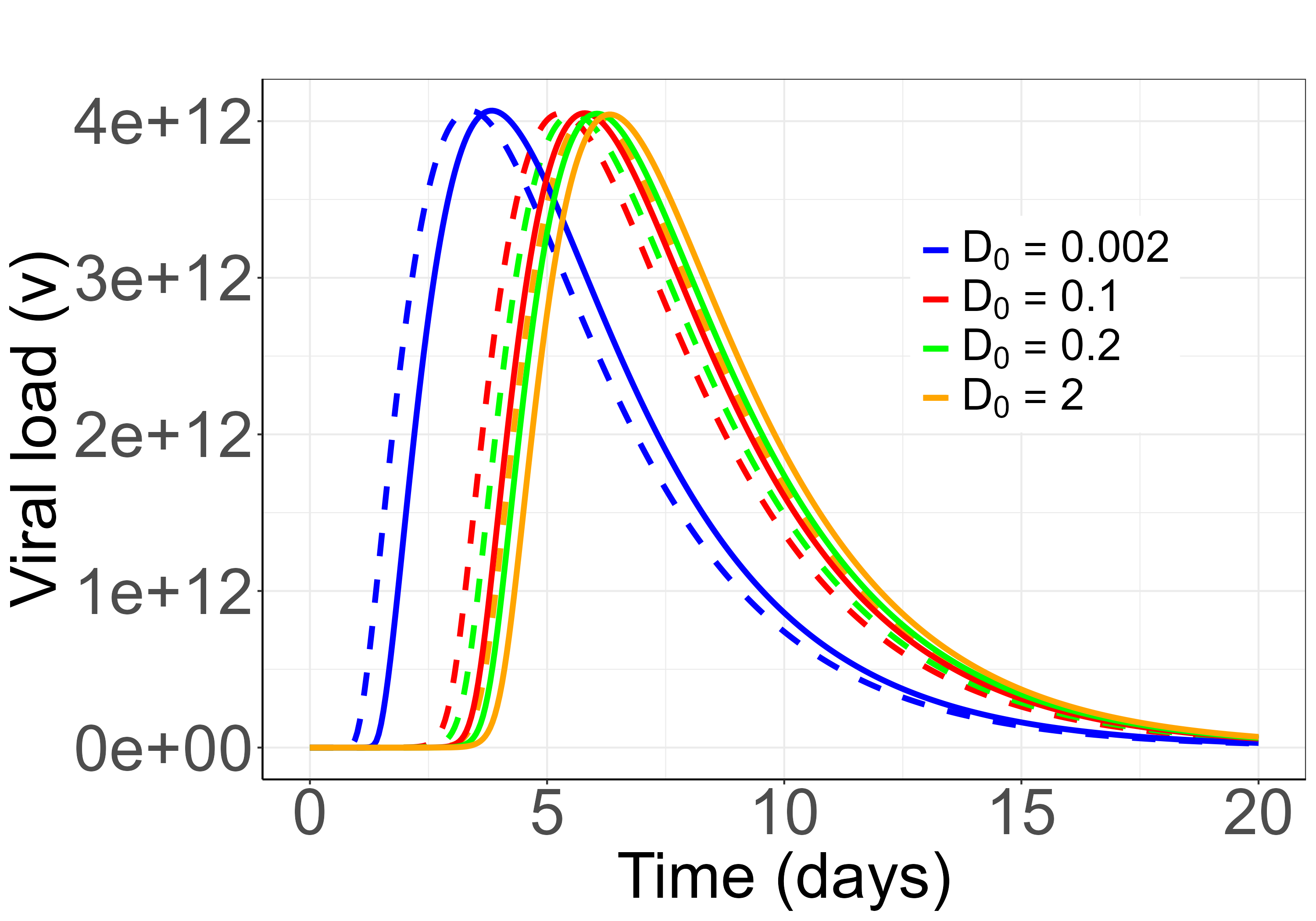}\\
  \includegraphics[width=2.0in]{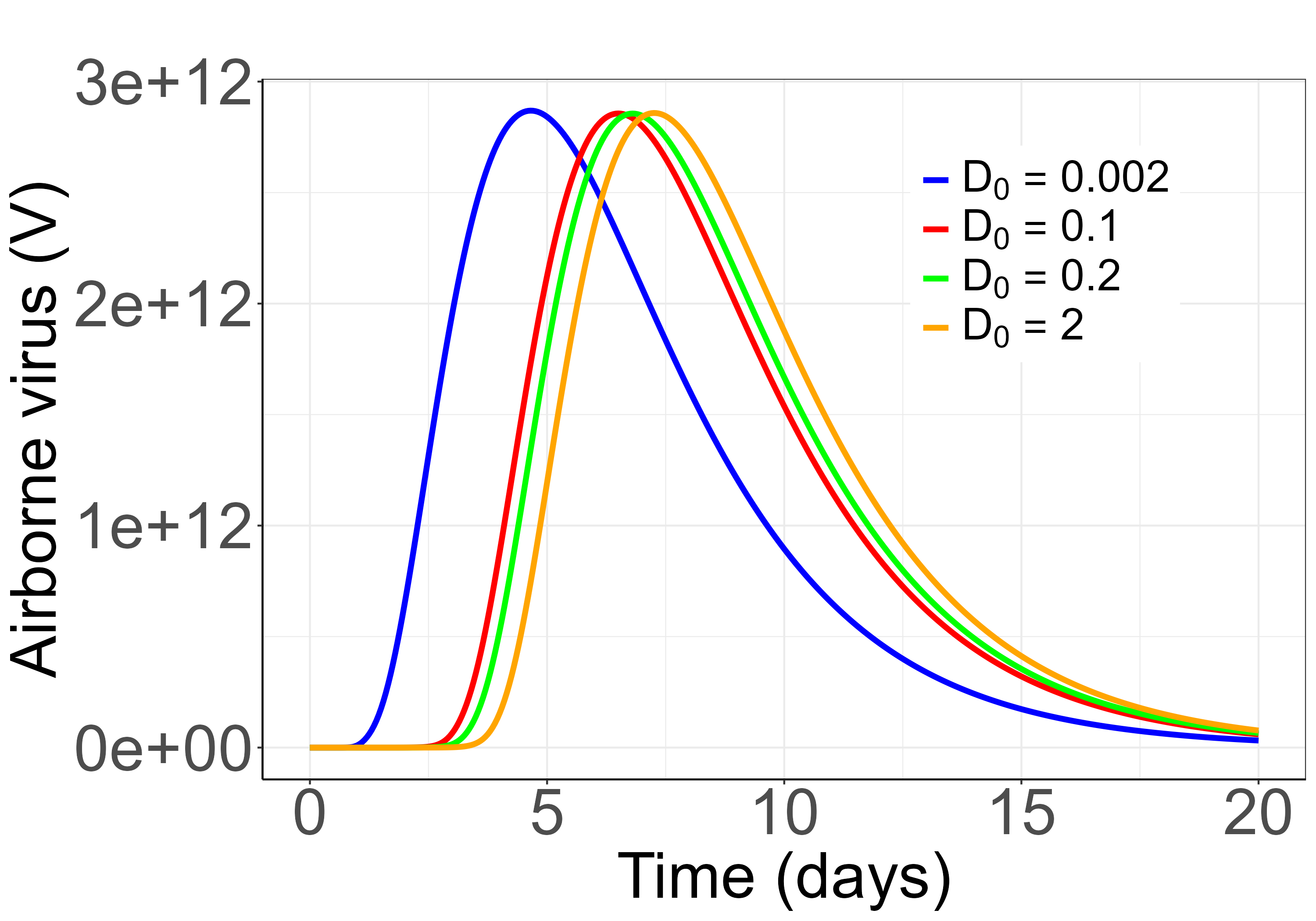} \quad
  \includegraphics[width=2.0in]{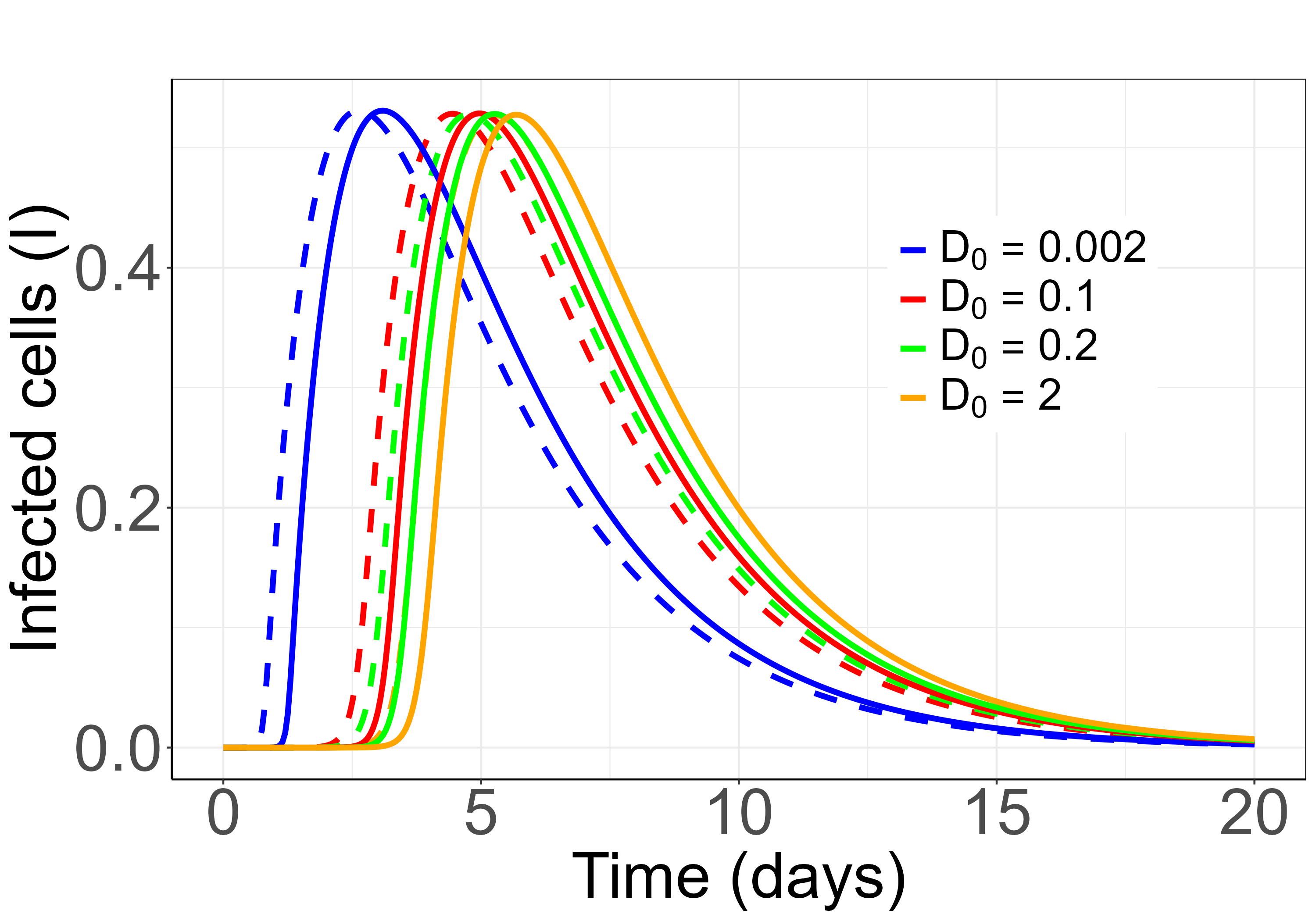} \quad
  \includegraphics[width=2.0in]{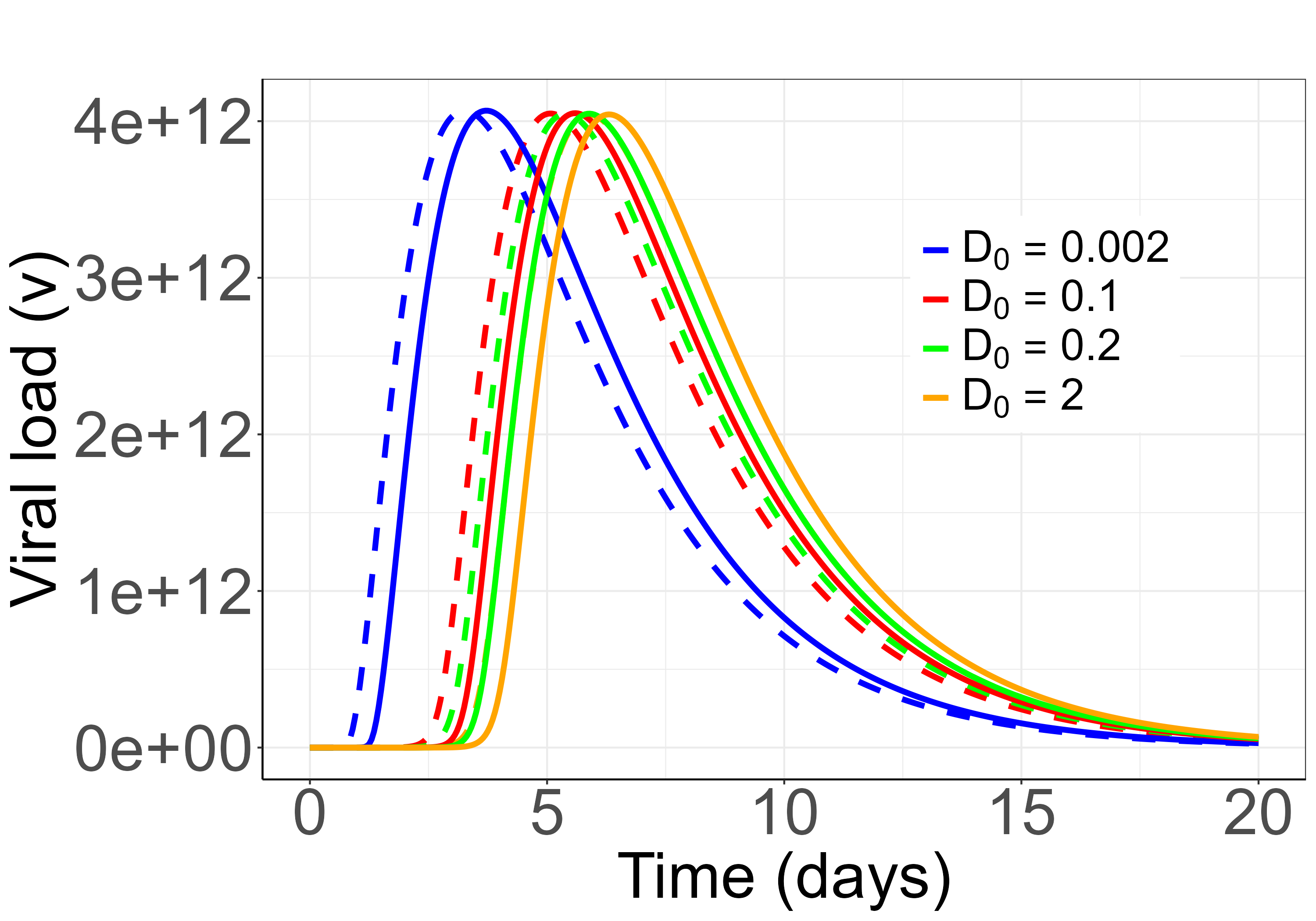}
  \caption{\textbf{Effects of diffusion and host position on within-host infection dynamics.}
Time series of airborne viral concentration (column~1), infected cells (column~2), and within-host viral load (column~3) computed using \eqref{eqn:TwoPopModel_PDE2}. The figure illustrates how diffusion and host position influence within-host infection dynamics for two identical hosts, $H_1$ (dashed lines) and $H_2$ (solid lines), with infection initiated in $H_1$. Host locations are shown by row: top row, $H_1 = (-0.15, 0)$ and $H_2 = (0.15, 0)$; middle row, $H_1 = (-0.5, 0)$ and $H_2 = (0.5, 0)$; and bottom row, $H_1 = (-0.85, 0)$ and $H_2 = (0.85, 0)$. All remaining parameters are listed in Table~\ref{Table:TableParameter}.}
\label{Fig:TwoHostDynamics}
\end{figure}

In this section, we consider a scenario involving multiple hosts within the domain. We investigate the impact of diffusing viral particles and host configuration on within-host infection dynamics. For simplicity, we focus on the case of two hosts in the domain and construct the corresponding multiscale model for this setting based on \eqref{eqn:ReducedODE_EclipDyn_IVP}, as described below.
\begin{equation}\label{eqn:TwoPopModel_PDE2}
\begin{split}
\frac{\text{d}V}{\text{d}t} &= \frac{1}{|\Omega|} (\xi_1 v_1 + \xi_2 v_2) - V , \\[2ex]
\frac{\text{d} T_1}{\text{d}t} &= -\beta_{11}\,T_1  \left( V + \frac{\xi_1 v_1}{2\pi D_0} \right)    - \beta_{21} T_1 v_1, \\
\frac{\text{d} E_{ 1}}{\text{d}t} & = \beta_{11}\,T_1  \left( V + \frac{\xi_1 v_1}{2\pi D_0} \right)  + \beta_{21}  T_1 v_1 - k_1 E_1,\\
\frac{\text{d} I_{1}}{\text{d}t} & = k_1 E_1 - \delta_1 I_1,\\
\frac{\text{d} {v}_1}{\text{d} t} & = \pi_1 {I}_1 - c_1 {v}_1 - \xi_1 \,v_1, 
\\[2ex]
\frac{\text{d} T_2}{\text{d}t} &= -\beta_{12}\,T_2  \left( V + \frac{\xi_2 v_2}{2\pi D_0} \right)  - \beta_{22} T_2 v_2, \\
\frac{\text{d} E_{ 2}}{\text{d}t} & =\beta_{12}\,T_2  \left( V + \frac{\xi_2 v_2}{2\pi D_0} \right)   + \beta_{22}  T_2 v_2 - k_2 E_2,\\
\frac{\text{d} I_{2}}{\text{d}t} & = k_2 E_2 - \delta_2 I_2,\\
\frac{\text{d} {v}_2}{\text{d} t} & = \pi_2 {I}_2 - c_2 {v}_2 - \xi_2 \,v_2.
\end{split}
\end{equation}
Using this model, we simulate infection kinetics in a two-host setting to examine how the diffusion rate of airborne viral particles and host separation influence within-host infection dynamics. We consider diffusion rates $D_0 = 0.001, 0.1, 0.2,$ and $2$, and analyze three spatial configurations for host locations: Case I, $H_1 = (-0.15, 0)$ and $H_2 = (0.15, 0)$, corresponding to an inter-host distance of $\Delta_H = 0.3$ units; Case II, $H_1 = (-0.5, 0)$ and $H_2 = (0.5, 0)$, with $\Delta_H = 1$ unit; and Case III, $H_1 = (-0.85, 0)$ and $H_2 = (0.85, 0)$, yielding $\Delta_H = 1.7$ units.
The results are presented in Figure~\ref{Fig:TwoHostDynamics}, which shows the dynamics of airborne viral particles (column 1), infected cells (column 2), and within-host viral load (column 3) for host~1 ($H_1$, dashed lines) and host~2 ($H_2$, solid lines) across different diffusion rates and host separations.
Throughout these simulations, the two hosts are assumed to be identical, with within-host parameter values taken from Table~\ref{Table:TableParameter}. Infection is initiated in host~1, with initial conditions $T_1(0)=1$, $E_1(0)=1/N(0)$, $I_1(0)=0$, and $v_1(0)=0$. Host~2 is initially uninfected, with $T_2(0)=1$, $E_2(0)=0$, $I_2(0)=0$, and $v_2(0)=0$. The initial airborne viral concentration is set to $V(0)=0$.

From Figure~\ref{Fig:TwoHostDynamics}, and consistent with the single-host results shown in Figure~\ref{Fig:SingleHostDynamics}, we observe a delay in the onset of within-host infection dynamics in both hosts. In the two-host setting, this delay is further influenced by the fact that infection is initiated in host~1, resulting in an additional lag before infection begins in host~2. This delay arises from the time required for host~1 to exhale viral particles and for these particles to diffuse through the air to reach host~2.
In Case~I, where the inter-host distance is smallest, infection in both hosts begins almost simultaneously (top row), particularly when the viral diffusion rate is low ($D_0 = 0.002$). However, as the diffusion rate increases, a noticeable delay in infection onset emerges in host~2. Moreover, as the inter-host distance increases (from the top to the bottom row), the delay in the onset of infection in host~2 becomes more pronounced.
Interestingly, while there is a clear difference in infection dynamics between Case~I (top row) and Case~II (middle row), the dynamics between Case~II and Case~III (bottom row) are largely similar, despite the increased inter-host distance. This suggests that beyond a certain separation threshold, further increases in distance have a diminishing impact on the timing of infection onset.

\section{Discussion}\label{sec:Discussion}

Integrating within-host infection dynamics with population-level disease transmission remains a significant challenge in mathematical epidemiology \cite{mideo2008linking,  sofonea2015within, almocera2018multiscale, feng2013mathematical, feng2012model, martcheva2015coupling}. In this article, we develop a novel multiscale modeling framework that bridges within-host viral replication with the spatial transport of virus-laden aerosols and population-level transmission processes. Our approach complements a growing body of work aimed at linking microscopic infection mechanisms to macroscopic outbreak dynamics, while providing a mechanistic and spatially explicit description of airborne transmission.
The proposed model describes within-host viral dynamics using a system of ODEs, coupled to a linear diffusion PDE that governs the spread of airborne viral particles in an enclosed environment. The coupling between scales is achieved through boundary conditions at the host-environment interface, representing viral shedding into, and inhalation from, the surrounding air. This formulation leads to a nonlocal coupled ODE-PDE system that explicitly connects individual-level infection kinetics with spatial transmission dynamics.

A key advantage of the proposed multiscale coupled ODE-PDE framework is its ability to capture the transient dynamics of both airborne viral particles and within-host infection processes for individual hosts. This allowed simultaneous characterization of transmission between hosts and the progression of infection within each host, providing a mechanistic understanding of disease spread across scales.
To enable analytical investigation, we applied matched asymptotic analysis to the coupled system and derived a reduced nonlinear system of ODEs in the regime of intermediate airborne viral diffusivity. This reduced model was more amenable to qualitative and quantitative analysis of the disease dynamics.
The derivation exploited the singular perturbation structure of the model through an asymptotic expansion in the small parameter $\varepsilon$, representing the ratio of host size to domain size. By constructing and matching inner and outer expansions around each host, the resulting multiscale ODE system retained essential spatial effects through the Neumann Green’s function and its regular part.

The derived multiscale ODE model eliminates the need to explicitly solve the diffusion equation at each time step while still retaining spatial heterogeneity. This formulation captures variations in the density of airborne viral particles across the domain and accounts for host spatial location. We establish the mathematical well-posedness of the model by proving the existence, uniqueness, and boundedness of solutions, thereby ensuring both its analytical validity and biological relevance.
Using this framework, we investigate within-host infection dynamics, beginning with the simplified setting of a single host in the domain. From the multiscale ODE model, we derive a two-term asymptotic expansion of the within-host basic reproduction number that explicitly incorporates spatial effects. This expansion reveals how domain geometry and host location influence infection potential. In the well-mixed limit ($D_0 \to \infty$), and when the contribution of diffusing airborne virus is neglected, the model reduces to the classical target cell limited model, recovering its standard expression for the basic reproduction number.
Numerical evaluation of the reproduction number further illustrates the role of diffusion in shaping within-host infection dynamics. When the diffusion rate is low, infection is primarily driven by inhalation of airborne viral particles concentrated near the host. As the diffusion rate increases, viral particles become more dispersed throughout the domain, and within-host viral dynamics increasingly dominate infection progression.

Furthermore, numerical simulation of the multiscale ODE model reveals that increasing the diffusion rate of airborne viral particles leads to a longer delay in the onset of infection within a host. This behavior arises because higher diffusivity disperses viral particles away from the immediate vicinity of the host, reducing the local concentration available for inhalation and subsequent infection. In the two-host scenario, where infection is initiated in one host, we observe that the time required for the second host to become infected is shorter when the inter-host distance is small, and increases as the distance between hosts grows. Notably, the delay in secondary infection also increases with the diffusion rate of airborne particles. Intuitively, one might expect that faster diffusion would enhance transmission by spreading particles more rapidly; however, our results suggest that a critical accumulation of airborne particles near a susceptible host is necessary to reach an infectious dose. High diffusivity prevents such accumulation, thereby delaying transmission. This finding is consistent with observations from the COVID-19 pandemic, where outdoor settings with greater air circulation and well-ventilated indoor environments have been associated with reduced risk of airborne transmission \cite{morawska2020can}.

Our multiscale modeling framework contributes to the growing literature on linking within-host infection kinetics with population-level transmission dynamics \cite{almocera2018multiscale, xue2020analysis, heitzman2022modeling} by explicitly incorporating the spatial dynamics of airborne virus transport in enclosed environments. The framework provides mechanistic and analytical insight into how indoor transmission is shaped by the diffusion of airborne viral particles and the spatial configuration of hosts. As such, it offers a flexible platform for evaluating integrated public-health and ventilation-based intervention strategies tailored to specific indoor settings, including classrooms, hospital wards, prison cells, and school dormitories.

Despite these strengths, our modeling framework has several limitations. First, it assumes that hosts remain fixed in space over the time scale of infection, whereas in reality individuals move within indoor environments as they carry out daily activities. Second, the model assumes that infected individuals continue their normal behavior and do not isolate after becoming infected, which may overestimate transmission in settings where isolation or behavioral changes occur. Finally, the diffusion of airborne viral particles is assumed to be spatially uniform throughout the domain, represented by a constant diffusion coefficient. In practice, airflow patterns, ventilation systems, and obstacles can induce spatially heterogeneous diffusion, which may significantly affect airborne transmission dynamics.

Future work should consider extending the model to incorporate spatially inhomogeneous diffusion rates, which would allow for more realistic representations of airflow patterns and ventilation effects in enclosed environments. While the analysis presented here focused primarily on understanding how diffusion and host spatial configuration influence within-host infection kinetics, an important next step is the development of methods for computing the basic reproduction number at the population level.
Given the novelty of the proposed multiscale framework, standard approaches for calculating the basic reproduction number may not be directly applicable, necessitating the development of alternative analytical or numerical techniques. Addressing this challenge would provide a critical link between individual-level infection dynamics and population-level outbreak potential.

Finally, we aim to apply this framework to investigate the transmission dynamics of SARS-CoV-2 in realistic experimental settings. A particularly relevant application involves controlled animal experiments in which hamsters are housed in adjacent cages, with one infected individual and neighboring animals monitored for subsequent infection \cite{sia2020pathogenesis, hawks2021infectious, heitzman2022modeling}. These studies have provided important insights into airborne transmission but have largely relied on ODE-based models that do not explicitly account for host spatial arrangement or the diffusion of virus-laden aerosols between cages \cite{heitzman2022modeling}. By incorporating spatial host locations and airborne viral transport, our multiscale framework offers a mechanistic extension of existing approaches and has the potential to yield new insights into the role of spatial structure and aerosol dynamics in SARS-CoV-2 transmission.

\section*{Funding}

The authors received no specific funding for this work.

\section*{Declaration of Competing Interest}

The authors declare no competing interests.


\bibliographystyle{unsrt}
\bibliography{References.bib}

\newpage
\renewcommand{\thesection}{S\arabic{section}}
\renewcommand{\thefigure}{S\arabic{figure}}
\renewcommand{\thetable}{S\arabic{table}}
\setcounter{section}{0}
\setcounter{figure}{0}
\appendixpageoff
\appendixtitleoff
\renewcommand{\appendixtocname}{Supplementary material}
	
\pagenumbering{arabic}
\renewcommand*{\thepage}{S-\arabic{page}}


\end{document}